\documentclass[11pt, reqno]{amsart}
\usepackage{amsfonts,amsmath,amsthm}
\usepackage{amssymb,epsfig}
\usepackage{pgf,tikz}
\usepackage{fancybox}
\usepackage[T1]{fontenc}
\usepackage{stmaryrd}
\usepackage{graphicx}
\usepackage{eso-pic}
\usepackage{fullpage}

\usepackage[unicode=true]{hyperref}
\hypersetup{colorlinks = true}
\hypersetup{
    colorlinks = true,
    linkcolor = {blue},
    citecolor={blue}
}

\usepackage[above,below,verbose,section]{placeins}
\makeatletter

\usepackage{mathrsfs}
\usepackage{mathpazo}
\usepackage{lastpage}
\usepackage[shortlabels]{enumitem}

\usepackage{tikz}
\usetikzlibrary{shapes.geometric, arrows, calc}
\tikzstyle{startstop} = [rectangle, rounded corners, minimum width=4cm, minimum height=2cm,text centered, draw=black, fill=red!30]
\tikzstyle{arrow} = [thick,->,>=stealth]
\makeatletter
\def\namedlabel#1#2{\begingroup
    #2%
    \def\@currentlabel{#2}%
    \phantomsection\label{#1}\endgroup
}
\makeatother

\newtheorem{thm}{Theorem}[section]

\newtheorem{lem}[thm]{Lemma}

\newtheorem{prop}[thm]{Proposition}

\newtheorem{rmk}[thm]{Remark}
\newcommand{\be}{\begin{eqnarray}}
\newcommand{\ee}{\end{eqnarray}}
\newcommand{\ben}{\begin{eqnarray*}}
\newcommand{\een}{\end{eqnarray*}}
\newcommand{\beal}{\begin{aligned}}
\newcommand{\enal}{\end{aligned}}
\newcommand{\beq}{\begin{equation}}
\newcommand{\eeq}{\end{equation}}
\numberwithin{equation}{section}

\newcommand{\eps}{\varepsilon}

\newcommand{\E}{\mathbb{E}}

\newcommand{\T}{\mathbb{T}}
\newcommand{\R}{\mathbb{R}}

\newcommand{\Z}{\mathbb{Z}}

\newcommand{\om}{\omega}
\newcommand{\dt}{\delta}
\newcommand{\Dt}{\Delta}

\newcommand{\cC}{\mathcal{C}}

\newcommand{\wt}{\widetilde }

\title{Convergence rate of periodic homogenization for Hamilton--Jacobi equations with supercritical viscosity}

\thanks{$\dagger$ {\it Statements and Declarations: }The authors declare no competing interests. No data is created or used in this work.}

\begin{document}

\maketitle

\centerline{\scshape Gengyu Liu$^*$,\; Zibo Wang$^\dagger$, \;Jianlu Zhang$^\ddagger$}
\medskip
{\footnotesize
    \centerline{ State Key Laboratory of Mathematical Sciences,}
    \centerline{Academy of Mathematics and Systems Science,}
    \centerline{Chinese Academy of Sciences, Beijing 100190, China}
    \centerline{{\it Email: }liugengyu@amss.ac.cn$^*$,\;zibowang@amss.ac.cn$^\dagger$,\; jellychung1987@gmail.com$^\ddagger$}
}
\medskip

\begin{abstract}
This paper is devoted to studying the convergence rate in periodic homogenization
of the viscous Hamilton--Jacobi equation (with supercritical viscosity index $\sigma>1$)
\[
\partial_t u_\eps+H(x/\eps,\partial_xu_\eps)=\eps^\sigma\Dt u_\eps,\quad  (x,t)\in\R^n\times\R_+,\
\]
subject to a given initial datum. For Tonelli Hamiltonians $H$, the solution  $u_\eps$ is proved to converge to a certain $u:\R^n\times\R\rightarrow\R$ (solving the effective problem) in the optimal rate $O(\eps^{\min\{\sigma-1,1\}})$. We also obtain a rate for the case of non-convex Hamiltonians and discussed its optimality therein. 
\end{abstract}

\section{Introduction}

For $T>0$, consider the viscous Hamilton--Jacobi equation 
\beq\label{eq:1}
\left\{
\begin{aligned}
&-\partial_t u_\eps+H(x/\eps,\partial_xu_\eps)=\eps^\sigma\Dt u_\eps,\\
&u_\eps(x,T)=\varphi(x),\quad (x,t)\in\R^n\times(-\infty,T], 
\end{aligned}
\right.
\eeq
with terminal datum $\varphi\in C^{0,1}(\R^n,\R)$ and Lipschitz continuous Hamiltonian $H:(x,p)\in \R^{2n}\rightarrow \R$ satisfying: 
\begin{itemize}
\item[\namedlabel{itm:C1}{$(\cC_1)$}] $H(\cdot,p)$ is $\Z^n-$periodic in $x$;
\item[\namedlabel{itm:C2}{$(\cC_2)$}] $H(x,\cdot)$ is coercive in $p$ for any $x$;
\item[\namedlabel{itm:C3}{$(\cC_3)$}] $\lim_{|p|\rightarrow+\infty}{\rm ess\inf}_{x\in\R^n}\{|H(x,p)|^2+(n+1)\partial_xH(x,p)\cdot p\}=+\infty$.
\end{itemize}
Evidence in \cite{LS, Car2, Tran-book} shows that, when $\sigma>1$, the viscosity solution $u_\eps$ of \eqref{eq:1} converges to $u\in C(\R^n\times\R_+,\R)$ locally uniformly as $\eps\rightarrow 0_+$, which solves 
the effective equation 
\begin{equation}\label{eq:hj-ef}
\left\{
\begin{aligned}
& -\partial_t u + \overline{H}(\partial_x u) = 0, \\
& u(x, T) = \varphi(x), \quad (x,t)\in\R^n\times(-\infty,T],
\end{aligned}
\right.
\end{equation}
with the effective Hamiltonian $\overline H: p\in \R^n\rightarrow\R$ identified as a special value that makes the associated {\it cell problem} solvable. Precisely, for any $p\in \R^n$, $\overline H(p)$ is the unique value making the following static equation 
\begin{equation}\label{eq:hj-c-lg}
H(x,p + d_x w) =  \overline{H}(p),\quad x\in \T^n(:=\R^n\slash\Z^n)
\end{equation}
solvable. The aim of this paper is to estimate the rate of convergence for the convergence of $u_\eps$ to $u$, i.e, an upper bound for $\|u_\eps-u\|_{L^\infty(\R^n\times(0,T])}$ as $\eps\rightarrow 0_+$. For Tonelli Hamiltonians, we are able to show that $\|u_\eps-u\|_{L^\infty(\R^n\times(0,T])}\leq O(\eps^{\min\{\sigma-1,1\}})$ by using the large deviation theory. Moreover, we construct examples to show the optimality of the rate. Even for non-convex Hamiltonians, we obtain the convergence rate that are optimal in some cases.

The quantitative homogenization was vastly studied in the past twenty years. For first-order static equations with multiscale homogenization, the convergence rate $O(\eps^{1/3})$ was firstly obtained in \cite{CD-Ishii}. By the same approach, \cite{Car1} obtained the rate of multiscale homogenization for static fully nonlinear elliptic equations with nonvanishing viscosity in \cite{Car2}, then discussed the different rates with different vanishing viscosity index $\sigma\in(0,+\infty)$.  For single-scale homogenization for equation 
\eqref{eq:1} with $\sigma=1$, \cite{QSTY} got $O(\eps^{1/2})-$rate and proved its optimality (under their context). Notice that the convexity of the Hamiltonian $H$ is not necessary in all above results, yet, an open problem was proposed in \cite{QSTY}, asking whether better rate in the case $\sigma=1$ can be obtained for general Hamiltonians with convexity.  A confirmed answer was recently given in \cite{LTY} by obtaining $O(\eps|\log\eps|)$ rate for a specific case with classical mechanical Hamiltonians (which allows the Hopf--Cole transformation). Our current consideration shares a similar interest with \cite{Car2}, and obtained an {\bf optimal} rate for the case $\sigma>1$. Technically, our conclusion is quite different from \cite{QSTY, LTY}, but mostly inherits from \cite{WZ} with necessary adaptation for the homogenization context.


It's remarkable that the optimal rate of convergence $O(\eps)$ for single-scale homogenization of convex first-order equations was recently obtained in \cite{TY}, and then $O(\eps^{1/2})$ for multi-scale case in \cite{HJ}. For earlier progress in this direction with nearly optimal rates of convergence, we refer the reader to \cite{Coop} and the references therein. In all these results, the convexity of Hamiltonians is essentially required, and an optimal control formula was used to get more quantitative properties of $\overline H(p)$. 
Different from all these results, we quantitatively considered a simultaneous vanishing-viscosity\slash homogenization process during the proof, and analysed their competition as $\eps\rightarrow 0_+$. As far as we know, such a treatment is also new. We also refer the authors to \cite{CD2,CG} for other consequences on the convergence rate of the vanishing viscosity process.

At last, we refer the authors to more quantitative results on the homogenization in different contexts, e.g., the quasi-periodic setting \cite{HTZ, LS}, the weakly coupled systems \cite{MT, MN} and PDEs with complicated boundary conditions \cite{HI, MN2}. 

\subsection{Main result}
We now describe our main results. Firstly, we should reinforce previous assumptions by the following:
\begin{itemize}
\item[\namedlabel{itm:C4}{$(\cC_4)$}] $H(x,p)\in C^2(T^*\T^n,\R)$;
\item[\namedlabel{itm:C5}{$(\cC_5)$}]$H(x,\cdot)$ is superlinear in $p$ for any $x\in \T^{n}$, i.e.
\[
\lim_{\vert p \vert \rightarrow \infty}H(x,p)/|p| = + \infty.
\] 
\item[\namedlabel{itm:C6}{$(\cC_6)$}] $H(x,\cdot)$ is positive definite in $p$ for any $x$;
\end{itemize}
Any Hamiltonian satisfying the above assumptions is called a {\it Tonelli type}. 

\begin{thm}\label{thm:1}

Assume \ref{itm:C1} - \ref{itm:C6} and $\varphi\in C^{0,1}(\R^n,\R)$. For $\eps\in(0,1)$ and fixed $T>0$, let $u_\eps$ be the viscosity solution to \eqref{eq:1} and $u$ be the viscosity solution to \eqref{eq:hj-ef} (with $\overline H$ defined in \eqref{eq:hj-c-lg}). 
There exists a constant $C_0:=C(H,n,\|\varphi\|_{C^{0,1}(\R^n)})>0$ such that 
 \beq
 \label{eq:thm1}
 \|u_\eps(x,t) -u(x,t)\|_{L^\infty(\R^n)} \le C_0(1+(T-t))\eps^{\min\{\sigma-1,1\}}.
 \eeq
 In particular, if $\varphi\in C^2(\R^n,\R)$ and $\|\varphi\|_{C^2(\R^n)}<\infty$, then there exists a constant $C_1:=C(H,n,\|\varphi\|_{C^2(\R^n)})>0$ such that
 \beq
  \|u_\eps(x,t) -u(x,t)\|_{L^\infty(\R^n)} \le  C_1(\min\{T-t,\eps\}+\eps^{\sigma-1}(T-t)).
 \eeq
\end{thm}
Meanwhile, the optimality of Theorem \ref{thm:1} can be verified by the following result.
\begin{thm}
For Hamiltonian $H(x,p)=\frac{1}{2}p^2+\cos x-1$ with $(x,p)\in\R\slash 2\pi\Z\times\R$, the associate homogenization \eqref{eq:1} with terminal datum $\varphi(x)\equiv0$ admits the following estimate 
\begin{equation}
\|u_\eps(x,0)-u(x,0)\|_{L^\infty(\R^n)}\ge c_1\eps^{\min\{\sigma-1,1\}},
\end{equation}
for some constant $c_1>0$. Here $u_\eps$ and $u$ are the viscosity solutions to \eqref{eq:1} and \eqref{eq:hj-ef}, respectively.
\end{thm}

\begin{rmk}
\begin{itemize}[itemsep=0.2em]

\item[(i)] In our recent preprint \cite{LWZ}, the optimal convergence rate of periodic homogenization for equation \eqref{eq:1} with subcritial viscosity $\sigma\in[0,1)$ was established by  double-variable approach and explicit examples. In that case, no convexity of Hamiltonians is required and lower regularity of Hamiltonians is sufficient. All these features benefit from the strong ellipticity of the assoicated $\overline H$, which is totally different from the definition in \eqref{eq:hj-c-lg} (see \cite{LWZ} for more details).

\item[(ii)]  For the critical case $\sigma=1$, the approach in this paper is no longer available. In the recent work \cite{LTY},  a spatially uniform rate $O(\eps|\log\eps|)$ was obtained for classical mechanical Hamiltonians. For the mechanical case, the Hopf--Cole transformation helps convert the Hamilton--Jacobi equation into a linear parabolic equation. Based on sharp kernel estimates obtained via the {\bf Doob transform}, a comparison is made between the resulting integral representation and the effective Hopf--Lax formula. The authors also raised the open question of whether the $O(\eps|\log\eps|)$ convergence rate remains valid for general strictly convex Hamiltonians or general periodic diffusions, which was recently claimed the solvability with the assistance of GPT-6 Astra in \cite{Qin}.

\item[(iii)] Quantitative homogenization of a multi-scale setting for $\sigma>1$ can be treated in the same spirit of Theorem \ref{thm:1}. On the other side, it is still open whether there exist strictly convex Hamiltonians, for which a lower bound for the associated multi-scale homogenization indicates the sharpness of the convergence rate; see \cite{HJ, LS} for relevant considerations. We leave these issues for future investigation.


\item[(iv)] The strict convexity and $C^2$ regularity of the Hamiltonians play an essential role for the optimal convergence rate. Evidence from \cite{CG} shows that the convergence rate of solutions in the vanishing viscosity process will be weakened if removing the strict convexity assumption. However, it remains open whether the $C^2$ regularity can be reduced to maintain the same optimal rate. 


\end{itemize}
\end{rmk}

For non-convex Hamiltonians, our approach can still be applied to get the convergence rate (may not be optimal) for the periodic homogenization in \eqref{eq:1}:
\begin{thm}[non-convex Hamiltonian]
\label{thm:3}
Assume $H\in C^{0,1}(T^*\T^n,\R)$ satisfying \ref{itm:C1}, \ref{itm:C2} and $\varphi\in C^{0,1}(\R^n,\R)$. For $\eps\in(0,1)$ and fixed $T>0$, let $u_\eps$ be the viscosity solution to \eqref{eq:1} and $u$ be the viscosity solution to \eqref{eq:hj-ef} (with $\overline H$ defined in \eqref{eq:hj-c-lg}). 
There exists a constant $C_2:=C(H,T,\|d_x\varphi\|_{L^\infty(\R^n)})>0$ such that
 \beq
 \label{eq:thm3}
  \|u_\eps(x,t) -u(x,t)\|_{L^\infty(\R^n\times[0,T])} \le \left\{
  \begin{aligned}
  &C_2\eps^{\min\{1/2,(\sigma-1)/2\}}, \quad & n=1, \\
  &C_2\eps^{\min\{1/3,(\sigma-1)/2\}}, \quad &n\ge 2.
  \end{aligned}
  \right.
 \eeq
 Such a rate in the case $n=1$ is optimal, in the sense that  there exist a Hamiltonian $H\in C^{0,1}(T^*\T,\R)$ of which 
\[
\|u_\eps-u\|_{L^\infty(\R^n\times[0,T])}\ge c_2\eps^{\min\{1/2,(\sigma-1)/2\}}
\]
 for some comstant $c_2>0$.
\end{thm}
Recently, \cite{JSTY} discussed the possible optimal convergence rate of periodic homogenization for non-convex Hamiltonians $H\in C^{0,1}(T^*\T^{n},\R)$ with $\sigma=+\infty$ (inviscid setting). For $n=1$, they improved  the $O(\eps^{1/3})$ rate in \cite{CD-Ishii} to $O(\eps^{1/2})$, with examples showing the sharpness. Moreover, they supply effective lower bounds of the rate: $O(\eps^{1/3}|\log\eps|^{-4/3})$ for $n=2$,  and $O(\eps^{1/3})$ for $n\ge3$. Inspired by their results, a similar argument as in our subsection \ref{s5} actually indicates the almost optimality of \eqref{eq:thm3} in the case $n=2$ and exact optimality in the case $n\geq 3$ both for $\sigma>5/3$. We don't know whether \eqref{eq:thm3} is sharp for $n\geq 2$ and $\sigma\in(1,5/3]$.

\medskip
\noindent{\bf Organization of the paper.} In Section \ref{s2}, we state some preliminary results on the viscous Hamilton--Jacobi equations and give a priori estimates for $u_\eps$ and $u$, which are used to simplify the settings of our problem. In Section \ref{s3} and Section \ref{s4}, we give the proof of Theorem \ref{thm:1} and Theorem \ref{thm:2}. In subsection \ref{s5}, we also show the convergence rate for non-convex Hamiltonians. We move a technical proof of regularity for solution to viscous Hamilton--Jacobi equations in Appendix \ref{appendix}.

\subsection*{Acknowledgements} 
This work is supported by the National Key R\&D Program of China (No. 2022YFA1007500) and the National Natural Science Foundation of China (No. 12231010, 11901560). ZW is additionally funded by the China Postdoctoral Science Foundation (No. 2024M763476) and supported by the National Center of Mathematics and Interdisciplinary Sciences. All the authors are grateful to Konstantin Khanin, Hung Vinh Tran and Hiroyoshi Mitake 
for helpful discussions and suggestions.

\section{Preliminary and setting}\label{s2}

\subsection{Gradient estimate}
For all $\sigma>1$, the effective Hamiltonian $\overline H:\R^n\rightarrow\R$ is continuous, convex and superlinear due to assumptions \ref{itm:C4}-\ref{itm:C6}. By the comparison principle, the solution $u(x,t)$ of \eqref{eq:hj-ef} satisfies 
\beq
\label{eq:bern}
\|\partial_t u\|_{L^\infty(\R^n\times(-\infty,T])}+\|\partial_x u\|_{L^\infty(\R^n\times(-\infty,T])}\leq M<+\infty
\eeq
for some constant $M:=M(\|d\varphi\|_{L^\infty(\R^n)},\overline H)>0$. Although in the above inequality it suffices to make $\varphi\in C^{0,1}(\R^n,\R)$, we indeed need more restrictive conditions to get a gradient estimate for $u_\eps$:

\begin{prop}[Bernstein method]
\label{prop-bernstein}
Suppose \ref{itm:C3}-\ref{itm:C6}, $\varphi\in C^2(\R^n,\R)$ and $\|\varphi\|_{C^2(\R^n)}<\infty$. There holds 
\beq\label{eq:bern-eps}
\|\partial_t u_\eps\|_{L^\infty(\R^n\times(-\infty,T])}+\|\partial_x u_\eps\|_{L^\infty(\R^n\times(-\infty,T])}\leq M'<+\infty
\eeq
for some constant $M':=M'(\|\varphi\|_{C^2(\R^n)}, H)>0$ independent of $\eps\in (0,1]$ and $\sigma\in[1,+\infty)$. 
\end{prop}

\subsection{Mollification argument for the terminal condition} 
Although the preceding gradient estimate \eqref{eq:bern-eps} is established under the terminal condition $\varphi\in C^{2}(\R^n,\R)$, the homogenization conclusion remains valid under a weaker regular terminal data $\varphi\in C^{0,1}(\R^n,\R)$. 
Indeed, the regularity assumption can be removed by a standard mollification argument (see, for instance, \cite{QSTY}). We include the proof for completeness.

Let $\rho\in C_c^\infty(\mathbb R^n)$ be a standard mollifier and
\[
\rho_\iota(x):=\iota^{-n}\rho(x/\iota), \quad
\varphi^\iota:=\rho_\iota * \varphi \quad \text{for }\iota>0.
\]
Then $\varphi^\iota\in C^2(\mathbb R^n,\R)$ and there exist constant $C_3:=C(\|d_x \varphi\|_{L^\infty(\R^n)})$ such that
\be
\label{est:moll}
\|\varphi^\iota-\varphi\|_{L^\infty(\R^n)}\le C_3\iota, \quad \|d_x\varphi^\iota\|_{L^\infty(\R^n)}\le C_3, \quad \|d_x^2\varphi^\iota\|_{L^\infty(\R^n)}\le \frac{C_3}{\iota}.
\ee
Let $u_\eps^\iota$ and $u^\iota$ be the solution to \eqref{eq:1} and \eqref{eq:hj-ef} with the modified terminal datum $\varphi^\iota\in C^2(\R^n,\R)$, respectively.
By \eqref{est:moll} and the comparison principle, 
we have
\be
\|u_\eps^\iota(x,t)-u_\eps(x,t)\|_{L^\infty(\R^n)}+\|u^\iota(x,t)-u(x,t)\|_{L^\infty(\R^n)}\le C_3\iota.
\ee
If the following Theorem \ref{thm:2} holds, we have
\be
\|u_\eps^\iota(x,t)-u^\iota(x,t)\|_{L^\infty(\R^n)}\le C(\min\{T-t,\eps\}+(T-t)\eps^{\sigma-1}).
\ee
Then we get the desired result by choosing $\iota=\eps^{\sigma-1}$ and 
\ben
\|u_\eps(x,t)-u(x,t)\|_{L^\infty(\R^n)}&\le& \|u_\eps(x,t)-u_\eps^\iota(x,t)\|_{L^\infty(\R^n)}\\
& &+\|u_\eps^\iota(x,t)-u^\iota(x,t)\|_{L^\infty(\R^n)}+\|u(x,t)-u^\iota(x,t)\|_{L^\infty(\R^n)},
\een
which finishes the assertion.

Therefore, to prove Theorem \ref{thm:1}, it is sufficient to show the following result. The proof will be given in the following Section \ref{s3}.
\begin{thm}
\label{thm:2}
Assume \ref{itm:C1} - \ref{itm:C6} and $\varphi\in C^{2}(\R^n,\R), \|\varphi\|_{C^2(\R^n)}<\infty$. Let $u_\eps$ be the viscosity solution to \eqref{eq:1} and $u$ be the viscosity solution to \eqref{eq:hj-ef} (with $\overline H$ defined in \eqref{eq:hj-c-lg}). There exists a constant  $C_1:=C(H,n,\|\varphi\|_{C^2(\R^n)})>0$ such that 
 \beq
  \|u_\eps(x,t) -u(x,t)\|_{L^\infty(\R^n)} \le C_1(\min\{T-t,\eps\}+\eps^{\sigma-1}(T-t)).
 \eeq
\end{thm}

\subsection{Optimal control formula} For any Hamiltonian $H(x,p)$ satisfying \ref{itm:C1}, \ref{itm:C5} and convexity in $p$,
the {\it Legendre transformation} establishes a Lagrangian function $L: (x,v)\in T\T^n\rightarrow \R$  by 
\beq\label{eq:led-tran}
L(x,v):=\max_{p\in\R^n}\big(\langle v,p\rangle-H(x,p)\big),
\eeq
which turns out to be continuous, convex and superlinear in $v$.  Moreover, if $H(x,p)$ satisfies \ref{itm:C4}-\ref{itm:C6}, the associated $L(x,v)$ can be verified to be  $C^2-$smooth and positive definite, superlinear in $v$ (also called {\it Tonelli Lagrangian}). 

The solutions to \eqref{eq:1} and \eqref{eq:hj-ef} possess the following optimal control interpretation:
\beq
\label{sto-op-con}
u_\eps(x,0)=\inf_{\substack{\gamma(0)=x,\\ d\gamma=\mathfrak v dt+\sqrt{2\eps^\sigma}dB_{\om}(t)}}\mathbb E\big\{\int_0^T L(\frac{\gamma}\eps,-\mathfrak v) dt+\varphi(\gamma(T))\big\},\quad \sigma\geq 1,
\eeq
and
\beq\label{det-op-con}
\quad u(x,0)=\inf_{\substack{\gamma(0)=x,\\ d\gamma=\mathfrak u dt}} \big\{\int_0^T \overline L(-\mathfrak u) dt+\varphi(\gamma(T))\big\}=\inf_{y\in\R^n} \{T\cdot \overline L((x-y)/T)+\varphi(y)\},
\eeq
where  $\mathfrak v\in L^1(\R^n,\R^n)$ could be any admissible progressively measurable control process, $\overline L(h)$  is the conjugated Lagrangian of $\overline H(p)$  as defined in \eqref{eq:led-tran} and the second equality in \eqref{det-op-con} is due to the Hopf--Lax formula. Here $\varphi\in C^{0,1}(\R^n)$ is sufficient to provide the representation.

To study the homogenization, it's helpful to involve a rescaling transformation by 
\[
(y,s):=(x/\eps,t/\eps)\in \R^{n+1}.
\]
Consequently, any curve $\gamma:t\in(-\infty,T]\rightarrow \R^n$ solving
\[
d\gamma (t)=\mathfrak vdt +\sqrt{2\eps^{\sigma}}dB_{\om}(t),\quad \gamma(0)=x
\]
will be transferred into a curve  $\eta:s\in(-\infty,T/\eps]\rightarrow \R^n$ solving 
\[
d\eta(s)=\mathfrak v ds+\sqrt{2\eps^{\sigma-1}}dB_{\om}(s), \quad \eta(0)=x/\eps,
\]
due to self-similarity of Brownian motion. 
Previous optimal control formula \eqref{sto-op-con} therefore becomes
\beq
\label{eq:stoc-op-eta}
u_\eps(x,0)=\inf_{\substack{\eta(0)=x/\eps,\\ d\eta=\mathfrak vds+\sqrt{2\eps^{\sigma-1}}dB_{\om}(s)}}\E \big\{\int_0^{T/\eps} \eps L(\eta, -\mathfrak v)ds+\varphi(\eps\eta(T/\eps))\big\}.
\eeq
This transformation allows a rescaled solution $w_\eps(y,s):=u_\eps(x,t)/\eps$ solving
\beq\label{eq:hj-vis-w}
\left\{
\begin{aligned}
&-\partial_s w_\eps+H(y,\partial_yw_\eps)=\eps^{\sigma-1}\Dt w_\eps,\\
&w_\eps(y,T/\eps)=\varphi(\eps y)/\eps,\quad (y,s)\in\R^n\times(-\infty,T/\eps].
\end{aligned}
\right.
\eeq
For simplify, we denote $\phi(y):=\varphi(\eps y)/\eps$ and $T':=T/\eps$.

\section{Convergence rate for convex case}\label{s3}

In this section, we give the proof of Theorem \ref{thm:2}.  The proof proceeds by separating the problem into two parts: the vanishing viscosity limit of the second-order Hamilton--Jacobi equation and homogenization of the first-order equation.
Recall that for the following equation:
\beq
\label{eq:hj-1st}
\left\{
\begin{aligned}
&-\partial_t \wt u_\eps+H(x/\eps, \partial_x\wt u_\eps)=0,\\
&\wt u_\eps(x,T)=\varphi(x),  \quad (x,t)\in\R^n\times(-\infty, T],
\end{aligned}
\right.
\eeq
its solution $\wt u_\eps$ will converge to the solution of \eqref{eq:hj-ef} with convergence rate of $O(\eps)$, see \cite[Theorem 1.1]{TY}. 
\begin{lem}[First-order homogenization]
\label{lem:der-hom}
Suppose that $H$ is convex in $p$ and satisfying \ref{itm:C1}-\ref{itm:C2}. Let $\wt u_\eps$ be the viscosity solution to \eqref{eq:hj-1st} and $u$ be the viscosity solution to \eqref{eq:hj-ef}. If $\varphi\in C^{0,1}(\R^n,\R)$, then there exists a constant $C_4:=C(\|d_x\varphi\|_{L^\infty(\R^n)})$ such that
\be
\|\wt u_\eps-u\|_{L^\infty(\R^{n})}\leq C_4\eps.
\ee
\end{lem}
By performing a scaling transformation similar to that in \eqref{eq:hj-vis-w}, we introduce the scaled function $\wt w_\eps(y,s)$ $:=\wt u_\eps(x,t)/\eps$ which satisfies
\beq\label{eq:hj-w}
\left\{
\begin{aligned}
&-\partial_s \wt w_\eps+H(y,\partial_y\wt w_\eps)=0,\\
&\wt w_\eps (y,T/\eps)=\varphi(\eps y)/\eps, \quad (y,s)\in\R^n\times(-\infty,T/\eps].
\end{aligned}
\right.
\eeq
Due to the classical weak KAM theory \cite{Ishii}, the solution to \eqref{eq:hj-w} is uniformly Lipschitz:
\be
\label{eq:lip-w}
\|\partial_s \wt w_\eps\|_{L^\infty(\R^n\times(-\infty,T/\eps])}+\|\partial_y \wt w_\eps\|_{L^\infty(\R^n\times(-\infty,T/\eps])}<M,
\ee
where we still denote the Lipschitz constant by $M$.

\subsection{Upper bound}
First, we use the uniformly semiconcavity of $w_\eps$ together with the comparison principle to obtain the upper bound.

\begin{lem}[Semiconcave]
\label{lem:semiconcave}
Let $w_\eps$ be the solution to \eqref{eq:hj-vis-w}. If $\varphi\in C^2(\R^n,\R)$ and $\|d_x^2\varphi\|_{L^\infty(\R^n)}<\infty$, then for every $s\in(-\infty, T']$, there exist a positive constant $C_\Dt:=C(\|d_x^2\varphi\|_{L^\infty(\R^n)},M',H)$ uniformly with respect to $\eps$ such that
\be
\label{semiconcave-w}
w_\eps(y+z,s)+w_\eps(y-z,s)-2w_\eps(y,s)\le C_\Dt |z|^2.
\ee
\end{lem}
\proof
For every $(y,s_0)\in\R^n\times(-\infty,T']$, we know that $w_\eps$ satisfies the stochastic optimal formula
\be
\label{eq:op-w-eps}
w_\eps(y,s_0)=\inf_{\substack{\eta(s_0)=y,\\ d\eta=\mathfrak vds+\sqrt{2\eps^{\sigma-1}}dB_{\om}(s)}}\E \big\{\int_{s_0}^{T'} L(\eta, -\mathfrak v)ds+\phi(\eta(T')\big\}.
\ee
Let $\mathfrak v_*(y,s):=-\partial_p H(y, \partial_y w_\eps(y,s))$ be the optimal vector field for which the infimum in \eqref{eq:op-w-eps} is achieved. 
By Proposition \ref{prop-bernstein}, the vector field $\mathfrak v_*$ is uniformly bounded, i.e., there exists a constant $C(M',H)$ such that $\|\mathfrak v_*\|_{L^\infty(\R^n\times(-\infty,T'])}<C(M',H)$.
When $s_0\le T'-1$, define a piecewise function
\begin{equation}
\chi(s)=\left\{
\begin{aligned}
&(s_0+1-s)z, &\quad s_0\le s\le s_0+1,\\
&0, &\quad s_0+1<s\le T',
\end{aligned}
\right.
\end{equation}
and 
\begin{equation}
\dot{\chi}(s)=\left\{
\begin{aligned}
&-z, &\quad s_0\le s\le s_0+1,\\
&0, &\quad s_0+1<s\le T'.
\end{aligned}
\right.
\end{equation}
For $s_0\le T'-1$, using the optimal control formula \eqref{eq:op-w-eps} and the optimality of $\mathfrak v_*$, we obtain
\be
\label{eq:semiconcave-1}
\begin{aligned}
&w_\eps(y+z,s_0)+w_\eps(y-z,s_0)-2w_\eps(y,s_0)\\
&\le \E\int_{s_0}^{s_0+1} L(\eta+\chi, -\mathfrak v_*-\dot{\chi})+L(\eta-\chi, -\mathfrak v_*+\dot{\chi}_1)-2L(\eta,-\mathfrak v_*)ds\\
&\le \E\int_{s_0}^{s_0+1} \|\partial_{xx}L\|\cdot |\chi|^2+2\|\partial_{xv}L\| \cdot|\chi||\dot{\chi}|+\|\partial_{vv}L\|\cdot|\dot{\chi}|^2 ds\\
&\le \sup_{|v|\le 1+\|\mathfrak v_*\|_{L^\infty}}(\|\partial_{xx}L\|+\|\partial_{xv}L\|+\|\partial_{vv}L\|) |z|^2\int_{s_0}^{s_0+1} (s_0+1-s)^2+(s_0+1-s)+1ds\\
&= \frac{11}{6}\sup_{|v|\le 1+\|\mathfrak v_*\|_{L^\infty}}(\|\partial_{xx}L\|+\|\partial_{xv}L\|+\|\partial_{vv}L\|) |z|^2.
\end{aligned}
\ee
Similarly, for $T'-1<s_0\le T'$, we have
\be
\label{eq:semiconcave-2}
\begin{aligned}
&w_\eps(y+z,s_0)+w_\eps(y-z,s_0)-2w_\eps(y,s_0)\\
&\le \E\int_{s_0}^{T'} L(\eta+z, -\mathfrak v_*)+L(\eta-z, -\mathfrak v_*)-2L(\eta,-\mathfrak v_*)ds\\
&+\E\big(\phi(\eta(T')+z)+\phi(\eta(T')-z)-2\phi(\eta(T'))\big)\\
&\le \E\int_{s_0}^{T'} \|\partial_{xx}L\|\cdot|z|^2ds+ \E\big(\|d_{x}^2\phi\|_{L^\infty(\R^n)}\cdot|z|^2\big)\\
&\le(\sup_{|v|\le1+\|\mathfrak v_*\|_{L^\infty}}\|\partial_{xx}L\|+ \eps\|d_x^2\varphi\|_{L^\infty(\R^n)})|z|^2.
\end{aligned}
\ee
Combining \eqref{eq:semiconcave-1} and \eqref{eq:semiconcave-2} together, we get the desired result \eqref{semiconcave-w}.
\qed

\begin{lem}[Upper bound]
\label{lem:upp-bound}
Let $w_\eps$ be the solution to \eqref{eq:hj-vis-w} and $\wt w_\eps$ be the solution to \eqref{eq:hj-w}. If $\varphi\in C^2(\R^n,\R)$ and $\|\varphi\|_{C^2(\R^n)}<\infty$, then for $t\in (-\infty, T]$
\be
\label{eq:upp-bound}
u_\eps(x,t)-\wt u_\eps(x,t) \le C_\Dt\eps^{\sigma-1}(T-t).
\ee
\end{lem}
\proof
Due to Lemma \ref{lem:semiconcave} we know that
\[
-\partial_s w_\eps+H(y,\partial_y w_\eps)\le \eps^{\sigma-1}C_\Dt.
\]
Therefore, $f(y,s):=w_\eps(y,s)-\eps^{\sigma-1}C_\Dt(T'-s)$ is a subsoution of \eqref{eq:hj-w}. By the comparison principle, this subsolution $f$ could not be bigger than $\wt w_\eps$ for every fixed $\eps>0$, i.e., 
\label{lem:upp-w}
\be
\label{eq:upp-bound}
w_\eps(y,s)-\wt w_\eps(y,s)\le C_\Dt \eps^{\sigma-1}(T'-s), \quad s\in(-\infty, T'].
\ee
Consequently, \eqref{eq:upp-bound} follows by $u_\eps=\eps w_\eps$.
\qed

\subsection{Lower bound}
While the upper bound relies on the uniform semiconcavity of the solution, the corresponding uniform semiconvexity estimate is available only on a finite initial time interval (Lemma \ref{lem:low-bound-hessian}). This allows us to obtain the lower bound of $w_\eps-\wt w _\eps$ on that interval by a similar argument as Lemma \ref{lem:upp-bound}. For longer times, the semiconvexity bound may deteriorate as $\varepsilon\to0$. We therefore introduce a convolution regular function and use an integrated estimate (Lemma \ref{lem:entropy-est}) to establish the lower bound in this regime.
\medskip

Due to the positive definiteness of $\partial_{pp}H$, there exists a constant $\theta>0$ such that 
\beq\label{eq:convexity}
\theta\mathbb{I}\preceq \partial_{pp}H(y,p),
\eeq
where $\mathbb I$ is the identity matrix of the rank $n$. Since $H\in C^2(T^*\T^n,\R)$, there exist constants $\Theta,\Theta_1,\Theta_2$ dependent on $M',H$ such that
 \be
 \label{eq:hessian-h}
 \|\partial_{pp}H\|_{L^\infty(\T^n\times B_0(M'))}\le \Theta, \quad  \|\partial_{yp}H\|_{L^\infty(\T^n\times B_0(M'))}\le \Theta_1, \quad  \|\partial_{yy}H\|_{L^\infty(\T^n\times B_0(M'))}\le \Theta_2.
 \ee

\begin{lem}
\label{lem:low-bound-hessian}
Let $w_\eps$ be the solution to \eqref{eq:hj-vis-w}. If $\varphi\in C^2(\R^n,\R)$ and $\|\varphi\|_{C^2(\R^n)}<\infty$, then there exist constants $s_*>0$ and $C_5:=C(\|d_x^2\varphi\|_{L^\infty(\R^n)})$ which are both uniform with respect to $\eps\in(0,1)$, such that for every $s\in(T'-s_*, T']$,
\be
\partial_{yy} w_\eps(s,y) \ge -C_5 \mathbb{I}.
\ee
Here $s_*>0$ is chosen strictly before the maximal existence time of the Riccati equation
\be
\label{eq:riccati}
k'(r)=\Theta k(r)^2+2\Theta_1 k(r)+\Theta_2, \quad k(0)=\eps\|d_x^2\varphi \|_{L^\infty(\R^n)}.
\ee
\end{lem}

\proof
For convenience, we first reversing the time $r=T'-s$ in \eqref{eq:hj-vis-w}:
\be
\label{eq:hj-vis-w-cov}
\left\{
\begin{aligned}
&\partial_r w_\eps+H(y, \partial_y w_\eps)=\eps^{\sigma-1}\Delta w_\eps, \\
& w_\eps(y,0)=\varphi(\eps y)/\eps, \quad (y,r)\in\R^n\times[0,\infty).
\end{aligned}
\right.
\ee
Differentiating equation \eqref{eq:hj-vis-w-cov} twice alone the eigenvector $\xi$, we get
\ben
\partial_r (\partial_{\xi\xi} w_\eps)+\langle \partial_p H(y,\partial_y w_\eps),\partial_y (\partial_{\xi\xi} w_\eps)\rangle+\partial_y(\partial_{\xi} w_\eps)^T\partial_{pp}H(y,\partial_{y} w_\eps)\partial_y(\partial_{\xi} w_\eps)\\
+2\xi^T\partial_{yp}H(y,\partial_y w_\eps)\partial_y(\partial_\xi w_\eps)+\xi^T\partial_{yy}H(y,\partial_y w_\eps)\xi=\eps^{\sigma-1}\Delta (\partial_{\xi\xi} w_\eps).
\een
Denote $m_\xi(r,y):=-\partial_{\xi\xi}w_\eps(r, y)$. Then using \eqref{eq:hessian-h} to obtain that
\be
\label{eq:hessian-est}
&&\partial_r m_\xi(r,y)+\langle \partial_p H(y,\partial_y w_\eps), \partial_y m_\xi(r,y)\rangle-\eps^{\sigma-1}\Delta m_\xi(r,y)\\
&&=\partial_y(\partial_{\xi} w_\eps)^T\partial_{pp}H(y,\partial_{y} w_\eps)\partial_y(\partial_{\xi} w_\eps)+2\xi^T\partial_{yp}H(y,\partial_y w_\eps)\partial_y(\partial_\xi w_\eps)+\xi^T\partial_{yy}H(y,\partial_y w_\eps)\xi \nonumber \\
&&\le \Theta |\partial_y(\partial_\xi w_\eps)|^2+2\Theta_1 |\partial_y(\partial_\xi w_\eps)|+\Theta_2 \nonumber.
\ee
In the following, we will use the maximum principle to find a uniform upper bound of $m_\xi$.
Denote 
\[
m(r):=\sup_{y\in\R^n, |\xi|=1} m_{\xi}(r,y).
\]
Note that the supremum of $m_\xi(r,y)$ may not be achieved on the non-compact space $\R^n$. Hence, we introduce a penalty function, i.e., for every $\delta>0$, let
\[
 m_{\delta}(r):=\sup_{y\in\R^n, |\xi|=1} \{m_\xi(r,y)-\delta\log(1+|y|^2)\},
\]
and denote the set of point achieved the maximum by $\mathcal{A}_\delta(r)$.
Since $\log(1+|y|^2)\to\infty$ as $|y|\to+\infty$, the set  $\mathcal{A}_\delta(r)$ is compact and $\lim_{\delta\to0}m_\delta(r)=m(r)$. Then, at each point $(y_\delta, \xi_\delta)\in \mathcal{A}_\delta(r)$, we have
\[
|\partial_y m_{\xi_\delta}(r,y_\delta)|= \frac{2\delta |y_\delta|}{1+|y_\delta|^2}\le \delta, \quad \Delta m_{\xi_\delta}(r,y_\delta)\le \delta\Big(\frac{2n}{1+|y_\delta|^2}-\frac{4|y_\delta|^2}{(1+|y_\delta|^2)^2}\Big)\le 2n\delta.
\]
By \eqref{eq:hessian-est}, we know that for every $r\ge 0$, there exists a constant $\Theta_3>0$ such that the following inequality holds
\be
\label{eq:der-m-xi}
\partial_r m_{\xi_\delta}(r, y_\delta)&\le&  \Theta |\partial_y(\partial_{\xi_\delta} w_\eps)|^2+2\Theta_1 |\partial_y(\partial_{\xi_\delta} w_\eps)|+\Theta_2+\Theta_3\delta\\
&=& \Theta|m_{\xi_\delta}(r,y_\delta)|^2+2\Theta_1|m_{\xi_\delta}(r,y_\delta)|+\Theta_2+\Theta_3\delta \nonumber,
\ee
where we used 
$|\partial_y(\partial_{\xi_\delta} w_\eps)|=|\partial_{yy}w_\eps \cdot \xi_\delta|=|m_{\xi_\delta}(r,y_\delta)|$.

By the regularity of $w_\eps$ (Lemma \ref{lem:regularity} in the Appendix), $\partial_r\partial_{yy}w_\eps$ is continuous for $r>0$. Hence the Danskin theorem \cite{Danskin} applies 
\[
\partial_+ m_{\delta}(r):=\lim_{h\to0_+}\frac{m_\delta(r+h)-m_\delta(r)}{h}=\max_{(y,\xi)\in\mathcal{A}_\delta(r)}\partial_rm_\xi(r,y).
\]
Set $M(r):=\max\{m(r),0\}$ and $M_\delta(r):=\max\{m_\delta(r),0\}$. If $m_\delta(r)\ge0$, then $0\le m_{\xi_\delta}(r,y_\delta)\le M(r)$ at every maximizing point. Taking positive parts and using \eqref{eq:der-m-xi}, we therefore obtain
\be
\label{eq:M_de}
\partial_+ M_\delta(r)\le \Theta M(r)^2+2\Theta_1M(r)+\Theta_2+\Theta_3\delta.
\ee
The regularity of $w_\eps$ also implies that 
\[
m_\delta,M_\delta\in C_{\mathrm{loc}}^{0,1}((0,\infty))\subset AC_{\mathrm{loc}}((0,\infty)).
\]
For any constant $a>0$, integrating the differential inequality \eqref{eq:M_de} on $[a,r]$ and letting $\delta\downarrow0$ gives
\[
M(r)\le M(a)+\int_a^r\bigl(\Theta M(q)^2+2\Theta_1M(q)+\Theta_2\bigr)\,dq.
\]

Finally, we explain the continuity of $M$ near $r=0$. The Duhamel formula gives that
\[
w_\eps(y,r)=G(y,r)*\phi-\int_0^r G(y,r-a)*H(y,\partial_y w_\eps(y,a))da,
\]
where $G$ is the heat kernel of the operator $\eps^{\sigma-1}\Delta$. Then by a direct calculation, we get
\[
\partial_{yy}w_\eps(y,a)\ge-\bigl(\eps\|d_x^2\varphi\|_{L^\infty(\R^n)}+C_\eps\sqrt a\bigr)I,
\]
uniformly in $y$ for sufficiently small $a>0$. Hence $\limsup_{a\downarrow0}M(a)\le\eps\|d_x^2\varphi\|_\infty$, and letting $a\downarrow0$ gives
\be
\label{eq:compar-M}
M(r)\le\eps\|d_x^2\varphi\|_{L^\infty(\R^n)}
+\int_0^r\bigl(\Theta M(q)^2+2\Theta_1M(q)+\Theta_2\bigr)\,dq.
\ee
Comparison \eqref{eq:compar-M} with \eqref{eq:riccati} now yields $m(r)\le M(r)\le k(r)$ before the maximal existence time of $k$. 
\qed

\begin{lem}[Lower bound in finite time interval]
\label{lem:short-lower}
Let $w_\eps$ be the solution to \eqref{eq:hj-vis-w} and $\wt w_\eps$ be the solution to \eqref{eq:hj-w}. If $\varphi\in C^2(\R^n,\R)$ and $\|\varphi\|_{C^2(\R^n)}<\infty$, then for $s\in(T'-s_*, T']$
\be
w_\eps(y,s)-\wt w_\eps(y,s)\ge -C_5\eps^{\sigma-1} (T'-s).
\ee
\end{lem}

\proof
Similar to Lemma \ref{lem:upp-w}, we know that $f(y,s):=w_\eps(y,s)+C_5 \eps^{\sigma-1}(T'-s)$ is a supersoultion of \eqref{eq:hj-w}, since
\[
-\partial_s w_\eps+H(y,\partial_y w_\eps)\ge -C_5\eps^{\sigma-1} .
\]
Using the comparison principle, we get the desired result.
\qed

\medskip

In the following, we introduce the sup-convolution of $\wt w_\eps$:
\beq
w^{\dt}(y,s):=\sup_{z\in\R^n} \Big(\wt w_\eps (z,s)-\frac{|z-y|^2}{\dt}\Big).
\eeq
\begin{prop}
\label{prop:sup-con}
The sup-convolution $w^\dt:\R^n\times(-\infty,T']\to\R$ satisfies the following properties
\begin{itemize}
\item[(1)] For every $s\le T'$, $0\leq w^{\dt}-\wt w_\eps\leq 2M^2\dt$;
\item[(2)] $w^\dt(\cdot, s)$ is $5M$-Lipschitz, $2C_{\Delta}$-semiconcave and $1/\delta$-semiconvex, i.e., $w^{\dt}+|y|^2/\dt$ is convex with respect to $y$;
\item[(3)]  There exist a positive constant $C_6:=C(M,H)$ such that $w^{\dt}$ is viscosity subsolution of
\beq
\label{eq:hj-w-dt}
\left\{
\begin{aligned}
-\partial_s w^{\dt}+H(y, \partial_y w^{\dt})&\leq C_6 \dt,\\
w^\dt(y,T')&=\phi_\dt(y) \leq \phi(y)+2\|d_y\phi\|_{L^\infty(\R^n)}^2\dt,
\end{aligned}
\right.
\eeq
where $\phi(y)=\varphi(\eps y)/\eps$ and $\|d_y\phi\|_{L^\infty(\R^n)}=\|d_x\varphi\|_{L^\infty(\R^n)}$.
\end{itemize}
\end{prop}
\proof
The proof is similar to that in \cite[Chapter 8]{Cal}, which treats the time-independent case, and \cite[Lemma 2]{Cecchin-Goffi}. We give the argument below for solutions of the Cauchy problem.\\
(1) For every $s\le T'$, if
\be
\label{eq:z-dt}
y^\dt\in{\arg\max}_{z\in\R^n}\Big(\wt w_\eps (z,s)-\frac{|z-y|^2}{\dt}\Big),
\ee
then 
\[
\wt w_\eps(y^\dt,s)-\wt w_\eps (y,s)=w^\delta(y,s)+\frac{|y^\dt-y|^2}{\dt}-\wt w_\eps (y,s)\ge\frac{|y-y^\dt|^2}{\dt}.
\]
Due to the Lipschitz property \eqref{eq:lip-w} of $\wt w_\eps$, we get $|y-y^\dt|\le M\dt$. For $(y,s)\in\R^n\times(-\infty,T']$,
\ben
|w^\dt(y,s)-\wt w_\eps (y,s)|&=&\Big|\wt w_\eps (y^\dt,s)-\frac{|y^\dt-y|^2}{\dt}-\wt w_\eps (y,s)\Big|\\
&\le& 2|\wt w_\eps(y^\dt,s)-\wt w_\eps(y,s)|\le 2M|y^\dt-y|\le 2M^2\dt.
\een
(2) For every $s\le T'$, if $|y-z|\ge M\dt$, then
\ben
|w^\dt(y,s)-w^\dt(z,s)|&\le& |w^\dt(y,s)-\wt w_\eps (y,s)|+|\wt w_\eps (y,s)-\wt w_\eps(z,s)|+|\wt w_\eps(z,s)-w^\dt(z,s)|\\
&\le&4M^2\dt+M|y-z|\le 5M|y-z|;
\een
if $|y-z|<M \dt$, let $y^\dt$ satisfying \eqref{eq:z-dt}, then
\ben
w^\dt(y,s)-w^\dt(z,s)&\le& \wt w_\eps(y^\dt,s)-\frac{1}{\dt} |y^\dt-y|^2-\wt w_\eps(y^\dt,s)+\frac1\dt|y^\dt-z|^2\\
&=&\frac{1}{\dt}(|y^\dt-z|^2-|y^\dt-y|^2)\le \frac 1\dt (|y-z|+2|y^\dt-y|)|y-z|\le 3M|y-z|.
\een
These two inequalities together give the Lipschitzness of $w^\delta$. Note that
\[
w^\delta+\frac{|y|^2}{\delta}=\sup_{z\in\R^n}\Big(\wt w_\eps(z,s)-\frac{|z|^2-2z\cdot y}{\delta}\Big),
\]
the convexity of $w^\delta+\frac{|y|^2}{\delta}$ follows by the property that the supremum of a family of affine functions is convex. The concavity of $w^\delta$ follows from the concavity of $\wt w_\eps$.\\
(3) Let $\psi(y,s)\in C^{2,1}(\R^n\times (-\infty,T'))$ be a test function. Suppose that $(\wt y,\wt s)$ is a local maximum point of $w^\dt-\psi$ and $z^\dt$ satisfies \eqref{eq:z-dt}. Hence
\be
\label{eq:w-dt-sub}
\wt w_\eps(z^\dt,\wt s)-\frac 1\dt|z^\dt-\wt y|^2-\psi(\wt y,\wt s)&=&w^\dt(\wt y, \wt s)-\psi(\wt y,\wt s)\\
&\ge& w^\dt (y, s)-\psi(y, s)\ge \wt w_\eps(z,s)-\frac1\dt|z-y|^2-\psi(y,s). \nonumber
\ee
Let $z=z^\dt$ we obtain that $\frac 1\dt|z^\dt-\wt y|^2+\psi(\wt y,\wt s)\le \frac1\dt|z^\dt -y|^2+\psi(y,s)$, which implies that $(\wt y,\wt s)$ is a local minimal point of $(y,s)\mapsto \frac1\dt |z^\dt-y|^2+\psi(y,s)$. Therefore, 
\be
\label{eq:psi_y}
|\partial_y\psi(\wt y,\wt s)|=\frac1\dt|\wt y-z^\dt|\le M.
\ee
Let $z^\dt-\wt y=z-y$ in \eqref{eq:w-dt-sub}, we have
\[
\wt w_\eps(z^\dt,\wt s)-\psi(\wt y,\wt s)\ge \wt w_\eps(z,s)-\psi(z-z^\dt+\wt y,s),
\]
which implies that $(z^\dt,\wt s)$ is a local maximum point of $(y,s)\mapsto \wt w_\eps(y,s)-\psi(z-z^\dt+\wt y,s)$.
Hence
\[
-\partial_s\psi(\wt y,\wt s)+H(z^\dt, \partial_y \psi(\wt y,\wt s))\le 0.
\]
By \eqref{eq:psi_y} and the regularity and periodicity of $H$, we have
\ben
-\partial_s\psi(\wt y,\wt s)+H(\wt y,\partial_y \psi(\wt y,\wt s))&\le&-H(z^\dt,\partial_y \psi(\wt y,\wt s))+H(\wt y,\partial_y \psi(\wt y,\wt s))\\
&\le& \|\partial_y H\|_{L^\infty(\T^n\times B_M(0))}\cdot|z^\dt-\wt y|\\
&\le& C(M,H)\dt,
\een
which shows that $w^\dt$ is a subsolution of \eqref{eq:hj-w-dt}.
\qed

\medskip

Subtracting \eqref{eq:hj-w-dt} from \eqref{eq:hj-w}, we obtain that
\beq
-\partial_s (w_\eps-w^\dt)+H(y,\partial_x w_\eps)-H(y,\partial_x w^\dt)\ge \eps^{\sigma-1}\Dt w_\eps- C_6\dt.
\eeq
By the fundamental theorem of calculus, the previous inequality becomes
\beq
\label{ineq-1}
-\partial_s (w_\eps-w^\dt)+\langle b(y,s), \partial_y (w_\eps-w^\dt)\rangle \ge \eps^{\sigma-1}\Dt w_\eps- C_6\dt,
\eeq
where the vector $b: \R^n\times(-\infty,T']\to\R^n$ is defined by
\be
b(y,s):=\int_0^1 \partial_p H(y, r\partial_y w_\eps(y,s)+(1-r)\partial_y w^\dt(y,s)) dr.
\ee
Let $\rho:\R^n\times(s_0,T']\to\R$ be the solution to the following Fokker--Planck equation:
\be
\label{eq:FP}
\left\{
\begin{aligned}
&\partial_s \rho(y,s)-\operatorname{div}(b(y,s)\rho(y,s))=\eps^{\sigma-1}\Dt\rho(y,s),\\
&\rho(y,s_0)=\dt_{y_0}, \quad  s_0< T'.
\end{aligned}
\right.
\ee

\begin{lem}
For every $(y,s)\in\R^n\times(s_0, T')$, there exist a constant $C_7:=C(n,\sigma, M,M')>0$ such that the solution to \eqref{eq:FP} satisfies
\be
\label{eq:entropy}
\Big|\int_{\R^n} \rho(y,s)\log\rho(y,s)  dy\Big|\le C_7(1+|\log\eps|+|\log s|).
\ee
\end{lem}
\proof
Using \eqref{eq:bern-eps} and Proposition \ref{prop:sup-con} 
we know that the vector $b$ is uniformly bounded, i.e.,
\[
\|b\|_{L^\infty(\R^n\times(-\infty,T'])}\le \|\partial_{p}H\|_{L^\infty(B(M'+5M)\times(-\infty,T'])}.
\]
By the upper bound estimate of $\rho$ (\cite{BKRS}, Corollary 7.2.3), for $\nu>(n+2)/2$ there exists a constant $C$ depending on $\nu,n$ such that 
\[
\rho(y,s)\le C(1+\eps^{1-\sigma})^\nu (1+s^{2\nu}\eps^{2\nu(1-\sigma)}) s^{-n/2}.
\]
Since $\rho$ is a probability density, we know that 
\[
\int_{\R^n} \rho(y,s) \log\rho(y,s)dy \le C(1+|\log s|+|\log\eps|),
\]
for $C:=C(\nu,n,\sigma)$.
To prove the lower bound, let $\overline \rho(y,s)$ be a Gaussian density which has the same second-order moment as $\rho(y,s)$, i.e.,
\[
\overline \rho(y,s):=\Big(\frac{n}{2\pi M_2(s)}\Big)^{n/2}\exp\Big(-\frac{n|y-y_0|^2}{2M_2(s)}\Big),
\]
where $M_2(s):=\int_{\R^n} |y-y_0|^2\rho(y,s)dy$. Due to the nonnegativity of relative entropy, we have
\be
\label{eq:entropy-1}
\int_{\R^n}\rho(y,s)\log\rho(y,s) dy&=& \int_{\R^n}\rho(y,s)\log\overline\rho(y,s) dy+ \int_{\R^n}\rho(y,s) \log\frac{\rho(y,s)}{\overline\rho(y,s)} dy\\ \nonumber
&\ge&\int_{\R^n}\rho(y,s)\log\overline\rho(y,s) dy \\ \nonumber
&=& \int_{\R^n} \Big(\frac{n}{2}\log n-\frac{n}{2}\log(2\pi M_2(s))-\frac{n|y-y_0|^2}{2M_2(s)}\Big)\rho(y,s)dy\\
&\ge& -\frac{n}{2}\log(2\pi M_2(s))+(\log n-\frac{n}{2}). \nonumber
\ee
To estimate the second-order moment $M_2$, we introduce a cutoff function to $f(y)=|y-y_0|^2$. For example, let $f_R(y)$ be a smooth function such that $f_R(y)=|y-y_0|^2$ for $|y-y_0|\le R$ and $f_R(y)=0$ for $|y-y_0|\ge 2R$. Moreover, there exists a constant $C_f>0$ independent on $R$ such that $|\partial_y f_R(y)|\le C_f f_R(y)^{1/2}$ and $\Dt f_R(y)\le C_f$. Then, applying \eqref{eq:FP} we have
\ben
\int_{\R^n} f_R(y)\rho(y,s)dy&=&\int_{s_0}^s\int_{\R^n} f_R(y) (\operatorname{div}(b(y,r)\rho(y,r))+\eps^{\sigma-1}\Dt\rho(y,r))dydr\\
&=&\int_{s_0}^s\int_{\R^n}-\langle\partial_y f_R(y), b(y,r)\rangle\rho(y,r)+\eps^{\sigma-1}\Dt f_R(y)\rho(y,r)dydr\\
&\le&\int_{s_0}^s \|b\|_{L^\infty(\R^n\times(-\infty,T'])} \Big(\int_{\R^n} C_f^2 f_R(y)\rho(y,r)dy\Big)^{1/2}+\eps^{\sigma-1}C_fdr\\
&=& \int_{s_0}^s \|b\|_{L^\infty(\R^n\times(-\infty,T'])}C_f\Big(\int_{\R^n} f_R(y)\rho(y,r)dy\Big)^{1/2}+\eps^{\sigma-1}C_fdr,
\een
where we used the integration by parts formula and H\"older inequality.
Letting $R\to\infty$ in the previous inequality, we obtain that
\[
M_2(s)\le \int_{s_0}^s  \|b\|_{L^\infty(\R^n\times(-\infty,T'])}C_f M_2(r)^{1/2}+\eps^{\sigma-1}C_f dr,
\]
since $\lim_{R\to\infty}f_R(y)=|y-y_0|^2$.
Then Gr\"ownall's inequality implies that
\[
M_2(s)\le 2\Big(\eps^{\sigma-1}C_f(s-s_0)+\|b\|_{L^\infty(\R^n\times(-\infty,T'])}^2C_f^2(s-s_0)^2\Big).
\]
Combining this together with the inequality \eqref{eq:entropy-1}, we know that there exists constant $C:=C(n,M,M')$ such that
\[
\int_{\R^n}\rho(y,s)\log\rho(y,s) dy\ge -C(1+|\log s|+|\log\eps|),
\]
which ends the proof.
\qed

\begin{lem}
\label{lem:entropy-est}
For every $\tau>0$, there exists a constant $C_8:=C(\theta,\Theta,n,M,M',C_\Dt)$ such that
\be
\int_{s_0+\tau}^{T'}\int_{\R^n} \Dt w^\dt(y,s)\rho(y,s)dyds \ge -C_8(1+|\log\eps|+|\log\tau|+|\log T'|+(T'-s_0)).
\ee
\end{lem}

\proof
Let $\chi_R(y):=\chi(y/R)$ be a nonnegative smooth cutoff function such that $\chi(r)=1$ for $|r|\le1$ and $\chi(r)=0$ for $|r|\ge2$.
Applying \eqref{eq:FP} and integration by parts, we have
\ben
&&\frac{d}{ds} \int_{\R^n} \log \rho(y,s)\rho(y,s) \chi_R(y) dy\\
&&=\int_{\R^n} \partial_s\rho(y,s)(\log\rho(y,s)+1) \chi_R(y)dy\\
&&=\int_{\R^n} \big(\operatorname{div}(b(y,s)\rho(y,s))+\eps^{\sigma-1}\Dt\rho(y,s)\big)(\log\rho(y,s)+1)\chi_R(y)dy\\
&&=\int_{\R^n} \langle -b(y,s), \partial_y\rho(y,s)\chi_R(y)\rangle + \langle -b(y,s)\rho(y,s), (\log\rho(y,s)+1)\partial_y \chi_R(y) \rangle dy\\
&&+\int_{\R^n} \eps^{\sigma-1} \Big(\Delta (\rho(y,s)(\log\rho(y,s)+1))-\frac{|\partial_y\rho(y,s)|^2}{\rho(y,s)}\Big)\chi_R(y)dy\\
&&=\int_{\R^n} \operatorname{div}(b(y,s))\rho(y,s)\chi_R(y) -\langle b(y,s) \rho(y,s), \log\rho(y,s)\partial_y\chi_R(y)\rangle dy\\
&&+\int_{\R^n} \eps^{\sigma-1} \rho(y,s)(\log\rho(y,s)+1)\Delta\chi_R(y)-\eps^{\sigma-1}\frac{|\partial_y\rho(y,s)|^2}{\rho(y,s)}\chi_R(y)dy,
\een
where we used $\Dt(\rho\log\rho)=\Dt\rho(\log\rho+1)+|\partial_y\rho|^2/\rho$ in the third equality.
Integrating the previous equality from $s_0+\tau$ to $T'$ we obtain that
\be
\label{eq:entropy-diss}
&&\int_{\R^n}\log \rho(y,T')\rho(y,T')\chi_R(y)dy-\int_{\R^n}\log \rho(y,s_0+\tau)\rho(y,s_0+\tau)\chi_R(y)dy\\
&&=\int_{s_0+\tau}^{T'}\int_{\R^n}\Big(\operatorname{div}(b(y,s))\rho(y,s) -\eps^{\sigma-1}\frac{|\partial_y\rho(y,s)|^2}{\rho(y,s)}\Big)\chi_R(y)dyds \nonumber \\ 
&&+\int_{s_0+\tau}^{T'}\int_{\R^n} -\langle b(y,s) \rho(y,s), \log\rho(y,s)\partial_y\chi_R(y)\rangle+\eps^{\sigma-1} \rho(y,s)(\log\rho(y,s)+1)\Delta\chi_R(y)dyds. \nonumber
\ee
Note that the cutoff function $\partial_y\chi_R=\partial_y\chi(y/R)/R$ and $\Delta\chi_R=\Delta\chi(y/R)/R^2$ are supported on $\{R<|y|<2R\}$, then we have
\ben
&&\Big|\int_{\R^n} -\langle b(y,s) \rho(y,s), \log\rho(y,s)\partial_y\chi_R(y)\rangle+\eps^{\sigma-1} \rho(y,s)(\log\rho(y,s)+1)\Delta\chi_R(y)dy\Big|\\
&&\le \Big|\int_{\{R<|y|<2R\}}\langle b(y,s) \rho(y,s), \log\rho(y,s)\partial_y\chi_R(y)\rangle dy\Big|\\
&&+\Big|\int_{\{R<|y|<2R\}}\eps^{\sigma-1} \rho(y,s)(\log\rho(y,s)+1)\Delta\chi_R(y)dy\Big|\\
&&\le \int_{\{R<|y|<2R\}} \Big( \frac{1}{R}\|b\|_\infty \|\partial_y\chi\|_\infty +\frac{1}{R^2} \|\Delta\chi\|_\infty\Big)|\rho(y,s)\log\rho(y,s)|dy,
\een
which converge to zero as $R\to\infty$, since the integral $\int_{\R^n}|\rho(y,s)\log\rho(y,s)|dy<\infty$ and  the vector field $b$ is bounded.
Therefore, letting $R\to\infty$ in the equality \eqref{eq:entropy-diss}, we get
\be
\label{eq:entr-ineq}
\begin{aligned}
\int_{\R^n}\log \rho(y,T')\rho(y,T')dy-\int_{\R^n}\log \rho(y,s_0+\tau)\rho(y,s_0+\tau)dy\\
\le \int_{s_0+\tau}^{T'}\int_{\R^n} \operatorname{div}(b(y,s))\rho(y,s)dyds.
\end{aligned}
\ee
From the definition of $b$ we get
\ben
\operatorname{div}(b(y,s))&=&\int_0^1 \big(\partial_y\cdot\partial_p\big)H(y,r\partial_y w_\eps+(1-r)\partial_y w^\dt)dr\\
&+&\int_0^1{\rm tr}\Big(\partial_{pp}H(y,r\partial_y w_\eps+(1-r)\partial_y w^\dt)\cdot\partial_{yy}(r w_\eps+(1-r) w^\dt)\Big)dr
\een
where $\partial_y\cdot\partial_p:=\sum_{i=1}^n\partial_{y_ip_i}$ for brevity.
By a property for trace of matrix (\cite{WZ}, Lemma 3.7), we have
\ben
&&{\rm tr}\Big(\partial_{pp}H(y,r\partial_y w_\eps+(1-r)\partial_y w^\dt)\cdot\partial_{yy}(r w_\eps+(1-r)w^\dt)\Big)\\
&&\le \theta {\rm tr}\partial_{yy}\big(r w_\eps+(1-r)w^\dt\big)+(\Theta-\theta)C_\Dt n=\theta (r\Dt w_\eps+(1-r)\Dt w^\dt)+(\Theta-\theta)C_\Dt n.
\een
Due to \eqref{eq:bern}, \eqref{eq:bern-eps} and the regularity of $H$, there exists a constant $C:=C(M,M',H,\Theta,\theta,C_\Delta,n)$ such that
\be
\label{eq:divb}
\operatorname{div}(b(y,s))\le C+\int_0^1 \theta (r\Dt w_\eps+(1-r)\Dt w^\dt) dr.
\ee
By \eqref{eq:entropy}, \eqref{eq:divb} and the semiconcavity of $w_\eps$, we get
\ben
\frac{\theta}{2}\int_{s_0+\tau}^{T'}\int_{\R^n} \Dt w^\dt(y,s)\rho(y,s)dyds &\ge& \int_{\R^n}\log \rho(y,T')\rho(y,T')dy-\int_{\R^n}\log \rho(y,s_0+\tau)\rho(y,s_0+\tau)dy\\
&&-\frac{\theta}{2}\int_{s_0+\tau}^{T'}\int_{\R^n} \Dt w_\eps(y,s)\rho(y,s)dyds-C(T'-s_0-\tau)\\
&\ge& -C_7(1+|\log\eps|+|\log T'|)-C_7(1+|\log\eps|+|\log(s_0+\tau)|)\\
&&-(C+\theta C_\Delta/2)(T'-s_0-\tau),
\een
which ends the proof.
\qed

\begin{lem}[Lower bound] 
\label{lem:low-w}
Let $w_\eps$ be the solution to \eqref{eq:hj-vis-w} and $w$ be the solution to \eqref{eq:hj-w}. If $\varphi\in C^2(\R^n,\R)$ and $\|\varphi\|_{C^2(\R^n)}<\infty$, then there exists a constant $C_9:=C(\theta,\Theta,n,M,M',C_\Dt,\|d_x\varphi\|_{L^\infty(\R^n)})$ such that
for any $t\in(-\infty,T]$,
\[
(u_\eps-\wt{u}_\eps)(x,t)\ge -C_9(1+(T-t))\eps^{\sigma-1}.
\]
\end{lem}
\proof
Define the smooth cutoff function $\chi_R(y):=\chi(y/R)$ as in Lemma \ref{lem:entropy-est}. Due to \eqref{ineq-1} and \eqref{eq:FP}, we have
\be
\label{eq:lower-1}
&&\int_{\R^n} (w_\eps-w^\dt)\rho(y,T')\chi_R(y)dy-\int_{\R^n}  (w_\eps-w^\dt)\rho(y,s_0)\chi_R(y)dy\\
&&=\int_{s_0}^{T'} \int_{\R^n} \partial_s(w_\eps-w^\dt)(y,s)\rho(y,s)\chi_R(y)+(w_\eps-w^\dt)(y,s)\partial_s\rho(y,s)\chi_R(y)dyds \nonumber\\
&&\le \int_{s_0}^{T'}\int_{\R^n} \Big(\langle b(y,s), \partial_y (w_\eps-w^\dt)\rangle- \eps^{\sigma-1}\Dt (w_\eps-w^\dt)- \eps^{\sigma-1}\Dt w^\dt+ C_3 \dt\Big)\rho(y,s)\chi_R(y)dyds \nonumber\\
&&+\int_{s_0}^{T'}\int_{\R^n} (w_\eps-w^\dt)\Big(\operatorname{div}(b(y,s)\rho(y,s))+\eps^{\sigma-1}\Dt\rho(y,s)\Big)\chi_R(y)dyds. \nonumber
\ee
Applying the integration by parts formula to the last integral in \eqref{eq:lower-1} yields that
\be
\label{eq:lower-2}
&&\int_{\R^n} (w_\eps-w^\dt)\Big(\operatorname{div}(b(y,s)\rho(y,s))+\eps^{\sigma-1}\Dt\rho(y,s)\Big)\chi_R(y)dy\\
&&=\int_{\R^n} -\Big(\langle \partial_y (w_\eps-w^\dt), b(y,s)\rangle \chi_R(y)+(w_\eps-w^\dt)\langle \partial_y\chi_R(y), b(y,s)\rangle\Big) \rho(y,s) dy \nonumber\\
&&+\int_{\R^n} \eps^{\sigma-1}\Big(\Delta(w_\eps-w^\dt) \chi_R(y)+(w_\eps-w^\dt) \Delta\chi_R(y)+2\langle \partial_y(w_\eps-w^\dt), \partial_y\chi_R(y)\rangle\Big)\rho(y,s)dy \nonumber\\
&&=\int_{\R^n} \Big(-\langle b(y,s),\partial_y (w_\eps-w^\dt)+\eps^{\sigma-1}\Delta(w_\eps-w^\dt)\Big) \rho(y,s)\chi_R(y)dy \nonumber\\
&&+\int_{\{R<|y|< 2R\}} \Big(-(w_\eps-w^\dt)\langle \partial_y\chi_R, b \rangle+(w_\eps-w^\dt)\Delta\chi_R+2\langle \partial_y(w_\eps-w^\dt), \partial_y\chi_R\rangle\Big)\rho(y,s)dy, \nonumber
\ee
where the second inequality holds since the cutoff function $\partial_y\chi_R, \Delta\chi_R$ are supported on $\{R<|y|<2R\}$. Note that the last integral in \eqref{eq:lower-2} converge to zero as $R\to\infty$ since $w_\eps-w^\dt, \partial_y(w_\eps-w^\dt)$ are bounded. Then combining \eqref{eq:lower-1}, \eqref{eq:lower-2} and letting $R\to\infty$, we get
\[
\int_{\R^n} (w_\eps-w^\dt)\rho(y,T')dy-\int_{\R^n}  (w_\eps-w^\dt)\rho(y,s_0)dy\le \int_{s_0}^{T'}\int_{\R^n} \Big(-\eps^{\sigma-1}\Delta w^\delta +C_6\delta\Big) \rho(y,s) dyds.
\]
For $s_0\le T'-\tau$, using the boundary conditions $\rho(y,s_0)=\dt_{y_0}$ and $(w_\eps-w^\dt)(y,T')=\phi(y)-\phi^\dt(y)$, we have
\ben
(w_\eps-w^\dt)(y_0,s_0)&\ge& 
\int_{\R^n} (\phi(y)-\phi_\dt(y))\rho(y,T')dy+\eps^{\sigma-1}\int_{s_0}^{s_0+\tau}\int_{\R^n}\Dt w^\dt(y,s)\rho(y,s)dyds\\
&+&\int_{s_0}^{s_0+\tau}\int_{\R^n} \Dt w^\dt(y,s)\rho(y,s)dyds-C_6\dt(T'-s_0).
\een
Due to Proposition \ref{prop:sup-con} and Lemma \ref{lem:entropy-est}, we get
\ben
(w_\eps-w^\dt)(y_0,s_0)&\ge&-2\|d_y\phi\|_{L^\infty(\R^n)}^2\dt -C_6(T'-s_0)\dt-2n \eps^{\sigma-1}\tau/\dt\\
&&-C_8\eps^{\sigma-1}(1+|\log\eps|+|\log\tau|+|\log T'|+(T'-s_0)).
\een
Choosing $\dt=\tau=\eps^{\sigma-1}$, we get for some constants $C_{10}>0$
\be
w_\eps-w^\dt(y_0,s_0)\ge 
\left\{
\begin{aligned}
&-C_{10}\eps^{\sigma-1}|\log\eps|, \qquad  \text{if } \eps^{\sigma-1}\le T'-s_0 \le |\log\eps|, \\
&-C_{10}\eps^{\sigma-1}(T'-s_0), \quad \text{if } T'-s_0 > |\log\eps|.
\end{aligned}
\right.
\ee
Combining this with Lemma \ref{lem:short-lower}, we get
\be
(w_\eps-\wt w_\eps)(y_0,s_0)\ge 
\left\{
\begin{aligned}
&-C_{5}\eps^{\sigma-1}(T'-s_0), \qquad \text{if } 0 \le T'-s_0< s_*,\\
&-C_{10}\eps^{\sigma-1}|\log\eps|, \qquad  \text{if } s_*< T'-s_0 \le |\log\eps|, \\
&-C_{10}\eps^{\sigma-1}(T'-s_0), \quad \text{if } T'-s_0 \ge |\log\eps|.
\end{aligned}
\right.
\ee
where we used $w^\dt-\wt w_\eps\ge0$.
The desired result follows by $u_\eps-\wt{u}_\eps=\eps(w_\eps-\wt w_\eps)$.
\qed

%

\medskip

\noindent {\it Proof of Theorem \ref{thm:2}.} 
Using the estimate from Lemma \ref{lem:der-hom}, Lemma \ref{lem:upp-bound} and Lemma \ref{lem:low-w}, we get
\[
u_\eps-u =(u_\eps-\wt{u}_\eps)+(\wt{u}_\eps-u) \leq C_\Dt\eps^{\sigma-1}(T-t)+C_4\eps,
\]
and 
\[
u_\eps-u =(u_\eps-\wt{u}_\eps)+(\wt{u}_\eps-u) \geq 
\left\{
\begin{aligned}
&-C_5\eps^{\sigma-1}(T-t)-C_4\eps, \qquad \text{if } 0 \le T-t\le \eps s_*,\\
&-C_{10}\eps^{\sigma}|\log\eps|-C_4\eps, \qquad  \text{if } \eps s_*\le T-t \le \eps|\log\eps|, \\
&-C_{10}\eps^{\sigma-1}(T-t)-C_4\eps, \qquad \text{if } T-t > \eps|\log\eps|.
\end{aligned}
\right.
\]
For the time approaching the terminal, using a strategy similar to \cite{HJ} yields an improved estimate. More precisely, applying the comparison principle we get 
\[
| u_\eps(x,t)-\varphi(x)|\le M' (T-t),
\]
and
\[
|u(x,t)-\varphi(x)|\le M (T-t).
\]
Consequently, for $t\in(-\infty,T]$,
\[
|u_\eps(x,t)-u(x,t)| \le |u_\eps(x,t)-\varphi(x)|+|u(x,t)-\varphi(x)|\le (M+M')(T-t).
\]
Therefore, the previous estimate can be improved to the following form
\be
|u_\eps-u|\le C(\min\{\eps,T-t\}+\eps^{\sigma-1}(T-t)),
\ee
for some constant $C>0$ uniformly with respect to $\eps$. This ends the proof.
\qed

\section{Examples and Optimality of the rates}\label{s4}

In this section, we show the optimality of the convergence rate by the example of a stochastic simple pendulum. More precisely, we consider the Hamiltonian $H(x,p)=\frac{1}{2}p^2+\cos x-1$ and $\sigma>1$.  Let $u_\eps$ and $u$ be the solution of 
\beq\label{eq:vis-hj-sp}
\left\{
\begin{aligned}
&-\partial_t u_\eps+\frac{1}{2}|\partial_xu_\eps|^2+\cos(x/\eps)-1=\eps^\sigma\Dt u_\eps,\\
&u_\eps(x,1)=0,\quad (x,t)\in\R\times(-\infty,1],
\end{aligned}
\right.
\eeq
and
\beq\label{eq:hj-ef-sp}
\left\{
\begin{aligned}
& -\partial_t u + \overline{H}(\partial_x u) = 0, \\
& u(x, 1) = 0, \quad (x,t)\in\R\times(-\infty,1],
\end{aligned}
\right.
\eeq
where $\overline H(p)$ is the unique value making the following static equation solvable
\begin{equation}
\frac{1}{2}|p + d_x w|^2+\cos x-1 =  \overline{H}(p),\quad x\in \R\slash 2\pi \Z.
\end{equation}
In the following, we sprarate the proof of Theorem \ref{thm:2} into two parts: $\sigma\in(1,2]$ and $\sigma\in(2,\infty)$.

\begin{lem}[{Optimality for $\sigma\in(1,2]$} ]
\label{lem:op-1}
Let $u_\eps$ and $u$ be solutions of \eqref{eq:vis-hj-sp} and \eqref{eq:hj-ef-sp} with the constant $\sigma\in(1,2]$. Then there exists a constant $c>0$ such that
\[
|u_\eps(0,0)-u(0,0)|\ge c \eps^{\sigma-1}.
\]
\end{lem}

\proof
According to \eqref{eq:hj-vis-w}, the rescale solution $w_\eps$ satisfies
\be
\left\{
\begin{aligned}
&-\partial_s w_\eps+\tfrac{1}{2}|\partial_y w_\eps|^2+\cos y-1=\eps^{\sigma-1}\Delta w_\eps, \\
&w_\eps(y,1/\eps)=0.
\end{aligned}
\right.
\ee
By the Hopf--Cole transformation, we know that $w_\eps(y,s)=-2\eps^{\sigma-1}\log \psi_\eps(y,s)$ and $\psi_\eps$ is the solution to the following backward Kolmogorov equation
\be
\left\{
\begin{aligned}
&\partial_s\psi_\eps+\tfrac{1}{2}\eps^{1-\sigma}(\cos y-1)\psi_\eps+\eps^{\sigma-1}\Delta\psi_\eps=0,\\
&\psi_\eps(y,1/\eps)=1.
\end{aligned}
\right.
\ee
Due to Feynman--Kac formula, the solution $\psi_\eps$ can be expressed as
\be
\label{eq:psi_eps}
\psi_\eps(0,0)=\E\Big(\exp\Big(\int_0^{1/\eps}\tfrac{1}{2}\eps^{1-\sigma}(\cos (X_r)-1) dr\Big)\Big), 
\ee
for $X_r=\sqrt{2\eps^{\sigma-1}}B_r$, in which $B_r$ is a standard one-dimensional Brownian motion.
Note that $-2\le \cos x-1\le 0$, we have $\psi_\eps\le 1$. In the following, we improve the upper bound of $\psi_\eps$.

Denote $\nu=\eps^{\sigma-1}$. First, we show that for $Y_r=y+\sqrt{2\nu}B_r$, there exists a constant $A\in(0,1)$ such that
\be
\label{eq:segment}
\sup_{y}\E_y\Big(\exp\Big(-\int_0^1 \tfrac{1}{2\nu}(1-\cos (Y_r)) dr\Big)\Big)\le A.
\ee
Let $\wt Y_t:=\inf_{k\in\Z}|Y_t-2k\pi|\in[0,\pi]$ be the distance of $Y_t$ to $2\pi\Z$, and let $\wt y:=\inf_{k\in\Z}|y-2k\pi|$. Then, for some constants $c_0>0$, we have
\be
\label{eq:est-segment}
&&\E\Big(\exp\Big(-\int_0^1 \tfrac{1}{2\nu}(1-\cos (Y_r)) dr\Big)\Big)\\
&=&\E\Big(\exp\Big(-\int_0^1 \tfrac{1}{2\nu}(1-\cos (\wt Y_r)) dr\Big)\Big) \nonumber\\
&
\le&\E\Big(\exp\Big(-\int_{0}^{1} \frac{|\wt Y_r|^2}{\nu\pi^2} dr\Big)\Big) \nonumber \\
&=&\E\Big(\Big(1_{\{\int_0^1|\wt Y_r|^2dr\ge c_0\nu \}}+1_{\{\int_0^1|\wt Y_r|^2dr< c_0\nu \}}\Big)\exp\Big(-\int_{0}^{1} \frac{|\wt Y_r|^2}{\nu\pi^2} dr\Big)\Big) \nonumber \\
&\le&\mathbb{P}\Big(\int_0^1|\wt Y_r|^2dr\ge c_0\nu\Big)\exp(-c_0/\pi^2)+\mathbb{P}\Big(\int_0^1|\wt Y_r|^2dr< c_0\nu\Big) \nonumber \\
&=&1-\mathbb{P}\Big(\int_0^1|\wt Y_r|^2dr\ge c_0\nu\Big)(1-\exp(-c_0/\pi^2)), \nonumber
\ee
where we used $1-\cos x\ge \frac{2}{\pi^2}|x|^2$ for $x\in[-\pi,\pi]$.
Indeed, we can show that there exist constants $c_0, p_0>0$ such that
\[
\mathbb{P}\Big(\int_0^1 |\wt Y_r|^2 dr\ge c_0\nu \Big)>p_0.
\]
We show this for two different cases: one with the initial point near $0$ and the other far from $0$: \\
(i) If the initial point $y$ satisfying $|\wt y|^2 \ge 2\nu$, then for $r\in[0,1]$
\ben
|\wt Y_r| = \inf_{k\in\Z} |y+\sqrt{2\nu}B_r-2k\pi| &\ge& \inf_{k\in\Z}| y-2k\pi| -\sqrt{2\nu}|B_r|\\
& \ge& |\wt y| -\sqrt{2\nu} \sup_{r\in[0,1]}|B_r| \ge \sqrt{2\nu}(1-\sup_{r\in[0,1]}|B_r|).
\een
By the reflection principle for Brownian motion, we know that
\ben
\mathbb{P}\Big( \sup_{r\in[0,1]} |B_r| \le  \sqrt{1/2}\Big)&\ge& 1-\mathbb{P}\Big(\sup_{r\in[0,1]} B_r >  \sqrt{1/2}\Big)-\mathbb{P}\Big(\inf_{r\in[0,1]} B_r <- \sqrt{1/2}\Big)\\
&=&1-2\mathbb{P}\Big(\sup_{r\in[0,1]} B_r >  \sqrt{1/2}\Big)\\
&=&1-4\mathbb{P}\Big(B_1 > \sqrt{1/2}\Big)=4\Phi( \sqrt{1/2})-3>0,
\een
where $\Phi$ is the distribution function for the standard normal distribution.
Therefore, 
\be
\label{eq:est-large-in}
\mathbb{P}_y\Big(\int_0^1 |\wt Y_r|^2 dr \ge(\sqrt{2}-1)^2\nu \Big) \ge \mathbb{P}\Big(\sup_{r\in[0,1]}|B_r|\le \sqrt{1/2}\Big)\ge 4\Phi( \sqrt{1/2})-3>0.
\ee

(ii) If the initial point $y$ satisfying $|\wt y|^2 < 2\nu $, then for $r\in[1/2,1]$
\ben
|\wt Y_r|= \inf_{k\in\Z} |Y_{1/2}+\sqrt{2\nu} (B_r-B_{1/2})-2k\pi|&\ge& \inf_{k\in\Z} |Y_{1/2}-2k\pi|-\sqrt{2\nu} |B_r-B_{1/2}|\\
&=&|\wt Y_{1/2}|-\sqrt{2\nu} |B_r-B_{1/2}|.
\een
Since the initial point $|\wt y|<\sqrt{2\nu} $, we know that when $2\le |B_{1/2}|\le 4$, the process $Y_{1/2}$ will not cross the periodic boundary. Hence, $|\wt Y_{1/2}|=|\wt y+\sqrt{2\nu}B_{1/2}|\ge \sqrt{2\nu}|B_{1/2}|-|\wt y|\ge \sqrt{2\nu }$ and
\ben
\mathbb{P}\Big(|\wt Y_{1/2}| \ge \sqrt{2\nu} \Big)&\ge& \mathbb{P}\Big(2 \le |B_{1/2}|\le 4\Big)=2(\Phi(4\sqrt{2})-\Phi(2\sqrt{2}))>0.
\een
Due to the Markov property and $B_r-B_{1/2}\stackrel{d}{=}B_{r-1/2}\stackrel{d}{=}B_{2r-1}/\sqrt{2}$, we have
\be
\label{eq:est-small-in}
\begin{aligned}
\mathbb{P}\Big(\int_0^1 |\wt Y_r|^2 dr \ge \nu /4\Big) &\ge \mathbb{P}\Big(|\wt Y_{1/2}| \ge \sqrt{2\nu} \Big)\mathbb{P}\Big(\sup_{r\in[1/2,1]}|B_r-B_{1/2}|\le 1/2\Big)\\
&=\mathbb{P}\Big(|\wt Y_{1/2}| \ge \sqrt{2\nu} \Big)\mathbb{P}\Big(\sup_{r\in[0,1]}|B_r|\le \sqrt{1/2}\Big)\\
&\ge 2(\Phi(4\sqrt{2})-\Phi(2\sqrt{2}))(4\Phi(\sqrt{1/2})-3)>0.
\end{aligned}
\ee
Combining \eqref{eq:est-segment}, \eqref{eq:est-large-in} and \eqref{eq:est-small-in}, we obtain that
\[
\E\Big(\exp\Big(-\int_0^1 \tfrac{1}{2}\eps^{1-\sigma}(1-\cos (Y_r)) dr\Big)\Big)\le 1-p_0 +p_0\exp(-c_0/\pi^2)<1,
\]
which shows the estimate \eqref{eq:segment}.

Using the equality \eqref{eq:psi_eps} and the Markov property of $X_r$, we have
\ben
\psi_\eps(0,0)&=&\E\Big(\exp\Big(\int_0^{1/\eps-1} \tfrac{1}{2}\eps^{1-\sigma}(\cos (X_r)-1) dr+\int_{1/\eps-1}^{1/\eps}\tfrac{1}{2}\eps^{1-\sigma}(\cos (X_r)-1) dr \Big)\Big)\\
&=&\E\Big(\exp\Big(\int_0^{1/\eps-1} \tfrac{1}{2}\eps^{1-\sigma}(\cos (X_r)-1) dr\Big)\\
&&\times\E\Big( \exp\Big(\int_{1/\eps-1}^{1/\eps}\tfrac{1}{2}\eps^{1-\sigma}(\cos (X_r)-1) dr\Big)\Big| X_{1/\eps-1}\Big)\Big)\\
&\le& A \E\Big(\exp\Big(\int_0^{1/\eps-1} \tfrac{1}{2}\eps^{1-\sigma}(\cos (X_r)-1) dr\Big)\Big).
\een
Consequently, $\psi_\eps(0,0)\le A^{\lfloor 1/\eps \rfloor}$
and 
\[
w_\eps(0,0)=-2\eps^{\sigma-1}\log\psi_\eps(0,0)\ge 2\eps^{\sigma-2}\log A.
\]
It is easy to get that $\overline{H}(0)=0$, and $u(0,0)=0$. Since $w_\eps(0,0)=-2\eps^{\sigma-1}\log\psi_\eps(0,0)\ge 0$, we have
\[
|u_\eps(0,0)-u(0,0)|=\eps w_\eps(0,0)-u(0,0)\ge 2\eps^{\sigma-1}\log A.
\]
This ends the proof.
\qed

\begin{lem}[Optimality for $\sigma\in(2,\infty)$ ]
Let $u_\eps$ and $u$ be solutions of \eqref{eq:vis-hj-sp} and \eqref{eq:hj-ef-sp} with the constant $\sigma\in(2,\infty)$. Then
\be
\label{eq:op}
|u_\eps(\pi\eps,0)-u(\pi\eps,0)|\ge 2 \eps.
\ee
\end{lem}

\begin{proof}
Consider the following first-order Hamilton--Jacobi equation
\beq
\left\{
\begin{aligned}
&-\partial_t \wt u_\eps+\frac{1}{2}|\partial_x\wt u_\eps|^2+\cos(x/\eps)-1=0,\\
&\wt u_\eps(x,1)=0,\quad (x,t)\in\R\times(-\infty,1].
\end{aligned}
\right.
\eeq
Its solution $\wt u_\eps$ has the following variational representation
\[
\wt u_\eps(x,0)=\inf_{\gamma\in AC([0,1]), \gamma(0)=x} \int_0^1 \Big(\frac12 |\dot{\gamma}(t)|^2+1-\cos(\gamma(t)/\eps)\Big) dt.
\]
First, there exists a zero-energy orbit 
\be
\dot{\gamma}(t)=-\sqrt{2(1-\cos(\gamma(t)/\eps))}, \quad \gamma(0)=\pi\eps,
\ee
such that
\ben
\int_0^1 \Big(\frac12 |\dot{\gamma}(t)|^2+1-\cos(\gamma(t)/\eps)\Big) dt&=&\int_0^1 2(1-\cos(\gamma(t)/\eps))dt\\
&=&\int_{\gamma(0)}^{\gamma(1)} -\sqrt{2(1-\cos( z /\eps))} dz\\
&=&\int_{\gamma(1)}^{\pi \eps}\sqrt{2(1-\cos( z /\eps))} dz\\
&\le& \int_{0}^{\pi \eps}\sqrt{2(1-\cos( z /\eps))} dz =4\eps.
\een
Since $\wt u_\eps$ is the infimum of the action, we have 
\be
\label{eq:wt-u-upp}
\wt u_\eps(\pi\eps, 0)\le 4\eps.
\ee
On the other hand, if $\gamma(t)$ is an orbit such that
\[
\int_0^1(1-\cos(\gamma(t)/\eps))dt\le  \int_0^1  \Big(\frac12 |\dot{\gamma}(t)|^2+1-\cos(\gamma(t)/\eps)\Big) dt <4\eps,
\] 
then there exists $t_0\in[0,1]$ such that
\be
\label{eq:orbit1}
1-\cos(\gamma(t_0)/\eps)\le 4\eps.
\ee
Therefore,
\be
\label{eq:action-low}
\begin{aligned}
\int_0^1 \Big(\frac12 |\dot{\gamma}(t)|^2+1-\cos(\gamma(t)/\eps)\Big) dt &\ge \int_0^{t_0} \Big(\frac12 |\dot{\gamma}(t)|^2+1-\cos(\gamma(t)/\eps)\Big) dt\\
&\ge \int_0^{t_0} |\dot{\gamma}(t)|\sqrt{2(1-\cos(\gamma(t)/\eps))} dt\\
&\ge \Big|\int_0^{t_0} \dot{\gamma}(t)\sqrt{2(1-\cos(\gamma(t)/\eps))} dt\Big|\\
&= |F(\gamma(t_0))-F(\gamma(0))|\\
&\ge |F(2k_0\pi\eps)-F(\pi\eps)|-|F(\gamma(t_0))-F(2k_0\pi\eps)|,
\end{aligned}
\ee
where $F(x):=\int_0^x \sqrt{2(1-\cos(z/\eps))}dz$ and $k_0\in\Z$ is chosen such that $2k_0\pi\eps$ is the closest potential well to $\gamma(t_0)$.
Since the interval of integration between $\pi\eps$ and $2k_0\pi\eps$ include at least half a period, we have
\be
\label{eq:F1}
\begin{aligned}
 |F(2k_0\pi\eps)-F(\pi\eps)|&=\int_{\pi\eps}^{2k_0\pi\eps}\sqrt{2(1-\cos(z/\eps))}dz \ge \int_0^{\pi\eps}\sqrt{2(1-\cos(z/\eps))}dz\\
 &=\eps\int_0^{\pi}\sqrt{2(1-\cos(z))}dz=4\eps.
 \end{aligned}
\ee
Using \eqref{eq:orbit1} we get 
\be
\label{eq:F2}
\begin{aligned}
|F(\gamma(t_0))-F(2k_0\pi\eps)|&=\int_{2k_0\pi\eps}^{\gamma(t_0)}\sqrt{2(1-\cos(z/\eps))}dz=\eps|\int_0^{\theta_0}\sqrt{2(1-\cos(z))}dz|\\
&=4\eps(1-\cos(|\theta_0/2|))\le 2\eps(1-\cos\theta_0)=2\eps(1-\cos(\gamma(t_0)/\eps))\le 8\eps^2,
\end{aligned}
\ee
where $\theta_0:=\gamma(t_0)/\eps-2k_0\pi$.
The inequalities \eqref{eq:F1} and \eqref{eq:F2} together with \eqref{eq:action-low} yields
\[
\int_0^1 \Big(\frac12 |\dot{\gamma}(t)|^2+1-\cos(\gamma(t)/\eps)\Big) dt\ge 4\eps-8\eps^2.
\]
By taking the infimum on the previous inequality, we get
\be
\label{eq:wt-u-low}
\wt u_\eps(\pi\eps, 0)\ge 4\eps-8\eps^2.
\ee
Combining \eqref{eq:wt-u-upp}, \eqref{eq:wt-u-low}, together with $u(x,t)=0$, we get
\[
|\wt u_\eps(\pi\eps,0)- u(\pi\eps,0)|=|\wt u_\eps(\pi\eps,0)|\ge 3\eps.
\]
Due to Lemma \ref{lem:upp-bound} and Lemma \ref{lem:low-w}, we know that for sufficiently small $\eps\in(0,1)$, 
\ben
|u_\eps(\pi\eps,0)-u(\pi\eps,0)|&\ge& |\wt u_\eps(\pi\eps,0)- u(\pi\eps,0)|-|u_\eps(\pi\eps,0)-\wt u_\eps(\pi\eps,0)|\\
&\ge& 3\eps-C_1\eps^{\sigma-1} \ge 2\eps,
\een
which ends the proof.
\end{proof}

 \subsection{Discussion of convergence rate for non-convex Hamiltonian}
 \label{s5}
 In this subsection, we discuss how the convergence rate changes when the strict convexity assumption on $H$ is dropped. The decomposition used in Section \ref{s3} still applies, i.e., it suffices to estimate the vanishing viscosity error and the first-order homogenization error separately.
\begin{lem}
\label{lem:non-con-1}
Assume \ref{itm:C1}-\ref{itm:C2} and $H\in C^{0,1}(\R^n\times\R^n)$. Let $\wt u_\eps$ be the viscosity solution to \eqref{eq:hj-1st} and $u$ be the viscosity solution to \eqref{eq:hj-ef}. If $\varphi\in C^{0,1}(\R^n,\R)$, then there exists a constant $C_{11}:=C(H,\|d_x\varphi\|_{L^\infty(\R^n)},T)$ such that
\be
\|\wt u_\eps-u\|_{L^\infty(\R^{n}\times[0,T])}\leq \left\{
\begin{aligned}
&C_{11}\eps^{1/2}, \quad n=1\\
&C_{11}\eps^{1/3}, \quad n\ge2.
\end{aligned}
\right.
\ee
\end{lem}

\begin{rmk}
The rate $O(\eps^{1/3})$ was obtained for arbritrary dimensions in \cite{CD-Ishii}. Recently,  for $n=1$, \cite{JSTY} improved the convergence rate to $O(\eps^{1/2})$ and verified the sharpness of the rate in the inviscid setting.

\end{rmk}


It remains to control the difference between the viscous and inviscid solutions.

\begin{lem}
\label{lem:non-con-2}
Assume \ref{itm:C1}-\ref{itm:C2} and $H\in C^{0,1}(\R^n\times\R^n)$. Let $u_\eps$ be the viscosity solution to \eqref{eq:1} and $\wt u_\eps$ be the viscosity solution to \eqref{eq:hj-1st}. If $\varphi\in C^{0,1}(\R^n,\R)$, then there exist a constant $C_{12}:=C(H,\|d_x\varphi\|_{L^\infty(\R^n)})$ such that
 \[
 \|u_\eps(x,t)-\wt u_\eps(x,t)\|_{L^\infty(\R^n)} \le C_{12}\eps^{(\sigma-1)/2}(T(\eps+T-t))^{1/2}.
 \]
\end{lem}
\begin{proof}
Let $w_\eps$ and $\wt w_\eps$ be solutions to \eqref{eq:hj-vis-w} and \eqref{eq:hj-w}. The vanishing viscosity estimate \cite[Theorem~3]{Cecchin-Goffi} implies that
 \be
 \|w_\eps(y,s)-\wt w_\eps(y,s)\|_{L^\infty(\R^{n})}\le C_{12}\eps^{(\sigma-1)/2}(T'(1+T'-s))^{1/2},
 \ee
 for some constant $C_{12}>0$.
Hence the desired result follows by $u_\eps-\wt u_\eps=\eps (w_\eps-\wt w_\eps)$.
\end{proof}

 \noindent {\it Proof of Theorem \ref{thm:3}.}  
 Due to Lemma \ref{lem:non-con-1} and Lemma \ref{lem:non-con-2}, there exist a constant $C_{13}$ depending on $H,\|d_x\varphi\|_{L^\infty(\R^n)},T,n$ such that
 \begin{equation*}
 \|u_\eps-u\|_{L^\infty(\R^n\times[0,T])}\le \|u_\eps-\wt u_\eps\|_{L^\infty(\R^n\times[0,T])}+ \|\wt u_\eps-u\|_{L^\infty(\R^n\times[0,T])}\le
 \left\{
  \begin{aligned}
  &C_{13}\eps^{\min\{1/2,(\sigma-1)/2\}}, \quad & n=1, \\
  &C_{13}\eps^{\min\{1/3,(\sigma-1)/2\}}, \quad &n\ge 2.
  \end{aligned}
  \right.
 \end{equation*}
We divide the proof of optimality at $n=1$ into two parts. The first part is dominant by the vanishing viscosity and the second part is dominant by first-order homogenization. 

(i) When $\sigma\in(1,2]$,  suppose that
\[
H(x,p):=(|p|-1)_+ -{\rm dist}(x,\Z),
\]
and $\varphi(x)=0$. Let $v:\T\to\R$ be a periodic solution to the following static equation
\be
\label{eq:static-sec}
H(y,v')=\overline{H}_\eps (0)+\eps^{\sigma-1} v''.
\ee
Then we can check that $U_\eps(x,t)=-\overline{H}_\eps (0)(T-t)+\eps v(x/\eps)$ satisfies $-\partial_t U_\eps+H(x/\eps,\partial_x U_\eps)=\eps^\sigma\Delta U_\eps$. Through a constant normalization $\inf_{x\in\T} v(x)=0$, we can check that $U_\eps^+:=U_\eps$ is a supersolution to \eqref{eq:1} and $U_\eps^-:=U_\eps-\eps\sup_{x\in\T} v$ is a subsolution to \eqref{eq:1}. Due to the comparision principle, 
\[
-\overline{H}_\eps(0)(T-t)+\eps (v(x/\eps)-\sup_{x\in\T}v)=U_\eps^-\le u_\eps\le U_\eps^+=-\overline{H}_\eps(0)(T-t)+\eps v(x/\eps)
\]
In the following, we show the lower bound of $-\overline{H}_\eps(0)$ by contradiction.
Denote $\chi(y):=v'(y)$. Due to \eqref{eq:static-sec} and its symmetric, we know that
\be
\label{eq:ex-nonconvex}
\eps^{\sigma-1}\chi'=(|\chi|-1)_+-y-\overline{H}_\eps(0), \quad 0<y<1/2,
\ee
and $\chi(0)=\chi(1/2)=0$.
If $-\overline{H}_\eps(0)\le \sqrt{2\eps^{\sigma-1}}$, then for sufficient small $\eps$, 
\[
\chi(y)=\eps^{1-\sigma} (-\overline{H}_\eps(0)y-y^2/2), \quad 0\le y\le -2\overline{H}_\eps(0)
\]
is a solution to \eqref{eq:ex-nonconvex}. This contradict with the fact that $\chi(y)>0$ on $y\in(0,1/2)$. Hence, $-\overline{H}_\eps(0)>\sqrt{2\eps^{\sigma-1}}$. Moreover, we obtain that
\[
(|\chi(y_*)|-1)_+=y_*+\overline{H}_\eps(0) <1/2, \quad 0<y_*<1/2,
\]
at the extreme point $y_*$ of $\chi$ (since $\chi'(y_*)=0$). This implies that $|\chi|<3/2$ and $\sup_{x\in \T} v \le 3/4$.

Recall that $\overline H(0)$ is the constant such that the follwing equation admits solution
\[
\overline{H}(0)=H(y,w').
\]
On the one hand, suppose that $y_0\in\T$ is the minimum point of $w$, i.e., $w(y_0)=\min_{y\in\T} w$.  Since the constant test function $\psi\equiv w(y_0)$ touches $w$ freom below at $y_0$, the viscosity supersolution property yields
\[
\overline{H}(0)\le H(y_0, 0)=-{\rm dist}(y_0,\Z)\le0.
\]
On the other hand, since $(|p|-1)_+\ge 0$, we have 
\[
\overline{H}(0)\ge -{\rm dist}(y,\Z), 
\]
for almost every $y\in\T$. Hence, $\overline{H}(0)=0$ and $u(x,t)\equiv 0$.
Therefore, for sufficient small $\eps\in(0,1)$, there exist a constant $C_{14}>0$ such that
\be
\label{eq:nonconvex-low-1}
u_\eps(x,0)-u(x,0)\ge -\overline{H}_\eps(0)\cdot T-\eps \sup_{x\in\T}v\ge C_{14}\eps^{(\sigma-1)/2}.
\ee

(ii) When $\sigma\in(2,+\infty)$, let $H$ be the one construct in \cite[Section 4.2]{JSTY} and $\varphi(x)=\min\{(x)_+,4\}$. Due to \cite[Theorem 4.4]{JSTY}, we know that 
\[
\liminf_{\eps\to 0}\frac{|\wt u_\eps(1,0)-u(1,0)|}{\eps^{1/2}}\ge \frac{1}{4\sqrt{2\pi}}.
\]
By Lemma \ref{lem:non-con-2}, we have $|u_\eps(1,0)-\wt u_\eps(1,0)|\le 2C_{12}\eps^{(\sigma-1)/2}$. Therefore, for sufficiently small $\eps\in(0,1)$, there exist a constant $C_{15}>0$ such that
\be
\label{eq:nonconvex-low-2}
\begin{aligned}
|u_\eps(1,0)-u(1,0)|&\ge |\wt u_\eps(1,0)-u(1,0)|-|u_\eps(1,0)-\wt u_\eps(1,0)|\\
&\ge  \frac{1}{16}\eps^{1/2}- 2C_{12}\eps^{(\sigma-1)/2}\ge C_{15}\eps^{1/2}.
\end{aligned}
\ee
Therefore, \eqref{eq:nonconvex-low-1} and \eqref{eq:nonconvex-low-2} gives the lower bound.
\qed

 \appendix
 \section{H\"older regularity of solutions for viscous Hamilton--Jacobi equation}
 \label{appendix}

In this Appendix, we show the regularity of the solution to the viscous Hamilton--Jacobi equation \eqref{eq:hj-vis-w-cov}. The argument is a standard regularity bootstrap. We include the details to verify that the stated assumptions suffice to ensure the continuity of the mixed derivative $\partial_{ryy}  w_\eps$. Here $C^{2k+\alpha, k+\alpha/2}$ is the classical H\"older space with respect to the parabolic distance. 
 \begin{lem}
 \label{lem:regularity}
Let $w_\eps$ be the solution to \eqref{eq:hj-vis-w-cov}. For every fixed $\eps, \alpha\in (0,1)$, $\sigma\in(1,\infty)$ and $0<\tau<T$,
\[
w_\eps, \partial_r w_\eps \in C^{2+\alpha,1+\alpha/2}(\R^n\times(\tau,T]).
\]
In particular, the classical mixed derivative exists for $r>0$
\[
\partial_{ryy} w_\eps\in C^{\alpha, \alpha/2}(\R^n\times(\tau,T]).
\]
\end{lem}

\begin{proof}
\noindent {\it Step 1: H\"older continuity of $\partial_y w_\eps$.} 
The Lipschitzness of $w_\eps$ gives that
\be
\|w_\eps\|_{L^\infty(\R^n\times[0,T])}\le \eps^{-1}\|\varphi\|_{L^\infty(\R^n)}+ M'T,
\ee
and
\be
\|\Delta w_\eps\|_{L^\infty(\R^n\times[0,T])}=\|\eps^{1-\sigma} (\partial_rw_\eps +H (y,\partial_y w_\eps(y,r)))\|_{L^\infty}\le \eps^{1-\sigma}(M'+C_{M'}).
\ee
For $1<p<\infty$, the interior elliptic $W^{2,p}$ estimate
\cite[Theorem~9.11]{GT} yields
\begin{equation}\label{eq:elliptic-1}
\begin{aligned}
\|w_\eps(\cdot,r)\|_{W^{2,p}(B_1(y_0))}
&\le C_{n,p}\bigl(\|w_\eps(\cdot,r)\|_{L^p(B_2(y_0))}
             +\|\Delta w_\eps(\cdot,r)\|_{L^p(B_2(y_0))}\bigr)\\
&\le C_{n,p}\left(\eps^{-1}\|\varphi\|_\infty+M'T+\eps^{1-\sigma}(M'+C_{M'})\right).
\end{aligned}
\end{equation}
Choosing $p>n/(1-\alpha)$, Morrey's inequality gives the spatial H\"older continuity of $\partial_y w_\eps$, i.e., 
\[
\sup_{0<r\le T}[\partial_y w_\eps(\cdot,r)]_{C^\alpha(\R^n)}\le C:=C(n,p,\eps,M',T).
\]
The time H\"older bound follows from this spatial estimate and time Lipschitzness. Indeed, for $0<\ell\le1$ and each coordinate
direction $e_i$, the forward spatial difference quotient gives
\[
|\partial_{y_i}w_\eps(y,r)-\partial_{y_i}w_\eps(y,r')|\le 2C\ell^\alpha+\frac{2M'|r-r'|}{\ell}.
\]
Taking $\ell=\sqrt{|r-r'|}$ when $0<|r-r'|\le1$ gives
\begin{equation}\label{eq:gradient-holder}
\|\partial_y w_\eps\|_{C^{\alpha,\alpha/2}(\R^n\times(0,T])}\le C.
\end{equation}
\medskip
\noindent {\it Step 2: Regularity improvement of $w_\eps$.} 
Let $f(y,r):=H(y,\partial_y w_\eps(y,r))$. Then
\ben
|f(y_1,r_1)-f(y_2,r_2)|&\le&|H(y_1,\partial_y w_\eps(y_1,r_1))-H(y_2,\partial_y w_\eps(y_1,r_1))|\\
&&+|H(y_2,\partial_y w_\eps(y_1,r_1))-H(y_2,\partial_y w_\eps(y_2,r_2))|\\
&\le& \|\partial_y H\|_{L^\infty(\T^n\times B_0(M'))} |y_1-y_2|\\
&&+\|\partial_p H\|_{L^\infty(\T^n\times B_0(M'))} |\partial_y w_\eps(y_1,r_1)-\partial_y w_\eps(y_2,r_2)|,
\een
which implies
\be
\label{eq:reg-f}
\|f\|_{C^{\alpha,\alpha/2}(\R^n\times(0,T])}\le C(H,M')\left(1+\|\partial_y w_\eps\|_{C^{\alpha,\alpha/2}(\R^n\times(0,T])}\right).
\ee
Then the Schauder estimate \cite[Theorem 8.11.1]{Krylov} improve the regularity
\begin{equation}\label{eq:schauder-1}
\|w_\eps\|_{C^{2+\alpha,1+\alpha/2}(\R^n\times[\tau/2,T])}
\le C\left(\|w_\eps\|_{L^\infty(\R^n\times(0,T])}+\|f\|_{C^{\alpha,\alpha/2}(\R^n\times(0,T])}\right)
\le C.
\end{equation}
\medskip
\noindent {\it Step 3: Regularity of $\partial_r w_\eps$ via time difference quotients.}
For fixed $0<\tau<T$, let $0<h<\tau/4$ and
\[
v_h(y,r):=\frac{w_\eps(y,r)-w_\eps(y,r-h)}{h}, \qquad (y,r)\in\R^n\times[\tau/2,T].
\]
Subtracting \eqref{eq:hj-vis-w-cov} at times $r$ and $r-h$ and using the mean value formula in the momentum variable yields
\begin{equation}\label{eq:app-time-difference}
\partial_rv_h+b_h(y,r)\cdot \partial_yv_h=\eps^{\sigma-1}\Delta v_h,
\end{equation}
where
\[
b_h(y,r):=\int_0^1 \partial_pH\bigl(y,
(1-\theta)\partial_yw_\eps(y,r-h)+\theta \partial_yw_\eps(y,r)\bigr)\,d\theta.
\]
The Lipschitzness of $w_\eps$ implies
\begin{equation}
\|v_h\|_{L^\infty(\R^n\times[\tau/2,T])}\le M', \quad |\partial_yw_\eps|\le M',
\end{equation}
which gives
\begin{equation}\label{eq:app-drift-holder}
\|b_h\|_{C^{\alpha,\alpha/2}(\R^n\times[\tau/2,T])}
\le C(H,M')\left(1+\|\partial_yw_\eps\|_{C^{\alpha,\alpha/2}(\R^n\times[\tau/4,T])}\right),
\end{equation}
uniformly for $0<h<\tau/4$.

For each fixed $h$, applying the Schauder estimate \cite[Theorem 8.11.1]{Krylov} to the linear parabolic equation \eqref{eq:app-time-difference} gives 
\begin{equation}\label{eq:app-schauder-time}
\|v_h\|_{C^{2+\alpha,1+\alpha/2}(\R^n\times[\tau,T])}
\le C\|v_h\|_{L^\infty(\R^n\times[\tau/2,T])}\le CM',
\end{equation}
where the constant $C$ depend on $\eps, \sigma, \tau, T, \alpha, n$ and the coefficient bounds in \eqref{eq:app-drift-holder}.
Since $\partial_rw_\eps$ is uniformly H\"older continuous in time, there exist constant $C$ independent on $h$ such that
\ben
|v_h(y,r)-\partial_rw_\eps(y,r)|&\le&\frac1h\int_{r-h}^r|\partial_rw_\eps(y,q) -\partial_r w_\eps(y,r)|dq \\
&\le& \frac Ch\int_0^h s^{\alpha/2} ds=\frac{C}{1+\alpha/2} h^{\alpha/2},
\een
for every $(y,r)\in \R^n\times[\tau,T]$. Therefore, $v_h$ uniformly converge to $\partial_r w_\eps$ and 
\[
\partial_rw_\eps\in C^{2+\alpha,1+\alpha/2}(\R^n\times[\tau,T]).
\]
Finally, on any positive-time strip, the identity
\[
w_\eps(y,r)-w_\eps(y,r_0)=\int_{r_0}^r \partial_q w_\eps (y,q)\,dq
\]
can be differentiated twice in $y$. Consequently,
\[
\partial_r\partial_y^2w_\eps=\partial_y^2\partial_rw_\eps\in C^{\alpha,\alpha/2}(\R^n\times[\tau,T]),
\]
which proves the asserted existence and continuity of the mixed derivative.
\end{proof}

\end{document}